\documentclass[11pt]{amsart}
\usepackage{amsmath}
\usepackage{amsxtra}
\usepackage{amscd}
\usepackage{amsthm}
\usepackage{amsfonts}
\usepackage{amssymb}
\usepackage{eucal}

\newtheorem{thm}{Theorem}[section]

\newtheorem{lem}[thm]{Lemma}
\newtheorem{prop}[thm]{Proposition}

\theoremstyle{remark}
\newtheorem{remark}[thm]{Remark}

\theoremstyle{definition}

\numberwithin{equation}{section}

\newcommand{\bean}{\begin{eqnarray}}
\newcommand{\eean}{\end{eqnarray}}
\newcommand{\be}{\begin{displaymath}}
\newcommand{\ee}{\end{displaymath}}
\newcommand{\bea}{\begin{eqnarray*}}   
\newcommand{\eea}{\end{eqnarray*}}

\newcommand{\thmref}[1]{Theorem~\ref{#1}}
\newcommand{\secref}[1]{Section~\ref{#1}}
\newcommand{\lemref}[1]{Lemma~\ref{#1}}
\newcommand{\propref}[1]{Proposition~\ref{#1}}

\newcommand{\remref}[1]{Remark~\ref{#1}}

\newcommand{\nc}{\newcommand}
\nc{\on}{\operatorname}
\nc{\ch}{\mbox{ch}}
\nc{\Z}{{\mathbb Z}}
\nc{\C}{{\mathbb C}}
\nc{\pone}{{\mathbb P}^1}
\nc{\pa}{\partial}
\nc{\F}{{\mathcal F}}
\nc{\arr}{\rightarrow}
\nc{\larr}{\longrightarrow}
\nc{\al}{\alpha}
\nc{\ri}{\rangle}
\nc{\lef}{\langle}
\nc{\W}{{\mathcal W}}
\nc{\la}{\lambda}
\nc{\ep}{\epsilon}
\nc{\su}{\widehat{{\mathfrak s}{\mathfrak l}}_2}
\nc{\sw}{{\mathfrak s}{\mathfrak l}}
\nc{\g}{{\mathfrak g}}
\nc{\h}{{\mathfrak h}}
\nc{\n}{{\mathfrak n}}
\nc{\N}{\widehat{\n}}
\nc{\G}{\widehat{\g}}
\nc{\De}{\Delta}
\nc{\gt}{\widetilde{\g}}
\nc{\Ga}{\Gamma}
\nc{\one}{{\mathbf 1}}
\nc{\z}{{\mathfrak Z}}
\nc{\La}{\Lambda}
\nc{\wt}{\widetilde}
\nc{\wh}{\widehat}
\nc{\cri}{_{\kappa_c}}
\nc{\kk}{h^\vee}
\nc{\sun}{\widehat{\sw}_N}
\nc{\si}{\sigma}
\nc{\el}{\ell}
\nc{\bi}{\bibitem}
\nc{\om}{\omega}
\nc{\ol}{\overline}
\nc{\ds}{\displaystyle}
\nc{\dzz}{\frac{dz}{z}}
\nc{\Res}{\on{Res}}
\nc{\mc}{\mathcal}
\nc{\Cal}{\mathcal}
\nc{\bb}{{\mathfrak b}}
\nc{\ot}{\otimes}
\nc{\R}{{\mc R}}
\nc{\yy}{{\mc Y}}
\nc{\ga}{\gamma}

\nc{\us}{\underset}
\nc{\opl}{\oplus}
\nc{\beq}{\begin{equation}}
\nc{\Fq}{{\mathcal F}}
\nc{\Mq}{{\mathcal M}}
\nc{\Rep}{\on{Rep}}
\nc{\sssec}{\subsubsection}
\nc{\ssec}{\subsection}
\nc{\lan}{\langle}
\nc{\ran}{\rangle}

\nc{\D}{\mathcal D}
\nc{\Vect}{\on{Vect}}
\nc{\ghat}{\G}
\nc{\T}{\mc T}
\nc{\Tloc}{\T^\g_{\on{loc}}}
\nc{\vac}{|0\ran}
\nc{\Wick}{{\mb :}}
\nc{\mb}{\mathbf}
\nc{\delz}{\partial_z}
\nc{\K}{{\cali K}}
\nc{\cali}{\mathcal}
\nc{\li}{\mathfrak l}
\nc{\lt}{\widetilde{\li}}
\nc{\astar}{a^*}
\nc{\cA}{{\mc A}}
\nc{\ka}{\kappa}

\nc{\OO}{{\mc O}}
\nc{\AutO}{\on{Aut}\OO}
\nc{\DerO}{\on{Der}\OO}
\nc{\DerpO}{\on{Der}_+\OO}
\nc{\Au}{{\mc A}ut}
\nc{\mf}{\mathfrak}
\nc{\V}{{\mathbb V}}
\nc{\hh}{\wh{\h}}

\nc{\pp}{{\mathfrak p}}
\nc{\mm}{{\mathfrak m}}
\nc{\rr}{{\mathfrak r}}
\nc{\ket}{\rangle}
\nc{\zz}{{\mathfrak z}}
\nc{\gr}{\on{gr}}
\nc{\Spe}{\on{Spec}}
\nc{\rv}{\crho}
\nc{\can}{\on{can}}
\nc{\CC}{\on{Op}_G(D))}
\nc{\Op}{\on{Op}_G(D)}
\nc{\MOp}{\on{MOp}_G(D)}
\nc{\Db}{{\mathbb D}}
\nc{\ww}{w}

\nc{\af}{{\mathbb A}^1}
\nc{\bs}{\backslash}
\nc{\laa}{(\la_i)}
\nc{\zn}{(z_i)}

\nc{\cla}{\check{\la}}
\nc{\cmu}{\check{\mu}}
\nc{\crho}{\check{\rho}}
\nc{\chal}{\check{\al}}
\nc{\cc}{{\mathfrak c}}

\nc{\tb}{\widetilde{\bb}}
\nc{\tn}{\widetilde{\n}}
\nc{\M}{{\mathbb M}}
\nc{\Nil}{{\mc N}}
\nc{\ppart}{(\!(t)\!)}

\begin{document}

\title{Self-extensions of Verma modules and differential forms on
  opers}

\author{Edward Frenkel}

\address{Department of Mathematics, University of California,
  Berkeley, CA 94720, USA}

\author{Constantin Teleman}

\address{DPMMS, Centre for Mathematical Sciences, Wilberforce Road,
  Cambridge, CB3 0WB, UK}

\date{January 2004, Revised: September 2005. MSC Primary 17B67,
  Secondary 17B56. Keywords: affine Kac-Moody algebra, critical level,
  Verma module, oper}

\begin{abstract}

We compute the algebras of self-extensions of the vacuum module and
the Verma modules over an affine Kac-Moody algebra $\ghat$ in suitable
categories of Harish-Chandra modules. We show that at the critical
level these algebras are isomorphic to the algebras of differential
forms on various spaces of opers associated to the Langlands dual Lie
algebra of $\g$, whereas away from the critical level they become
trivial. These results rely on and generalize the description of the
corresponding algebras of endomorphisms obtained by Feigin and Frenkel
and the description of the corresponding graded versions due to
Fishel, Grojnowski and Teleman.

\end{abstract}

\maketitle


\section*{Introduction}

Let $\g$ be a simple finite-dimensional Lie algebra and $\ghat_{\ka}$,
where $\ka$ is an invariant inner product on $\g$, the corresponding
affine Kac-Moody algebra. Consider the vacuum module $\V_\ka$ over
$\ghat_{\ka}$ (see \secref{vac module} for the precise
definitions). According to the results of \cite{FF:gd,F:wak}, the
algebra of endomorphisms of $\V_\ka$ is trivial, i.e., isomorphic to
$\C$, unless $\ka=\ka_c$, the critical value. In contrast, the algebra
$\on{End}_{\ghat_{\ka_c}} \V_{\ka_c}$ is large and is in fact
canonically isomorphic to the algebra of functions on the space
$\on{Op}_{^L \g}(D)$ of $^L \g$--opers on the disc, where $^L \g$ is
the Lie algebra that is Langlands dual to $\g$ (its Cartan matrix is
the transpose of that of $\g$).

In this paper we consider the algebra of endomorphisms of $\V_{\ka}$
in the derived category of $(\ghat_{\ka},G[[t]])$--modules, or in
other words, the algebra of self-extensions of $\V_{\ka}$ in the
abelian category of $(\ghat_{\ka},G[[t]])$--modules (here $G$ is the
connected simply-connected algebraic group corresponding to $\g$). As
we show in \secref{vac module}, this algebra may be realized as the
relative cohomology
\begin{equation}    \label{rel coh}
H^\bullet(\g\ppart,\g,\on{End}_\C \V_{\ka}) \simeq
H^\bullet(\g[[t]],\g,\V_{\ka}).
\end{equation}
We show that for $\ka \neq \ka_c$ we have $H^i(\g\ppart,\g,\on{End}_\C
\V_{\ka}) = 0$ if $i>0$, so the corresponding algebra of
self-extensions of $\V_{\ka}$ is isomorphic to $\C$. But if
$\ka=\ka_c$, then this algebra is isomorphic to the algebra of
differential forms on $\on{Op}_{^L \g}(D)$ (more precisely, this
isomorphism is defined up to a scalar which is fixed once we choose an
invariant bilinear form $\ka_0$ on $\g$). Moreover, deforming
$\V_{\ka_c}$ to $\V_\ka$, we obtain a differential on the cohomology
$H^\bullet(\g\ppart,\g,\on{End}_\C \V_{\ka_c})$ which coincides with
the de Rham differential.

In order to prove this result, we use the quasi-classical statement
about the cohomology $H^\bullet(\g[[t]],\g,\V^{\on{cl}})$, where
$\V^{\on{cl}}$ is the associate graded of $\V_{\ka_c}$ with respect to
the PBW filtration. The space $\on{Op}_{^L \g}(D)$ is an affine space
modeled on the vector space $C_{\g^*,\om}$ defined in \secref{qc
case}. According to \cite{FGT}, $H^\bullet(\g[[t]],\g,\V^{\on{cl}})$
is isomorphic to the algebra of differential forms on
$C_{\g^*,\om}$. Using this result and the description of
$H^0(\g[[t]],\g,\V_{\ka_c})$ from \cite{FF:gd,F:wak} mentioned above,
we obtain our result by employing the spectral sequence on the
cohomology corresponding to the PBW filtration.

We also prove an analogue of this result for the Verma modules
$\M_{\la,\ka}$. In this case we consider the algebra of endomorphisms
of $\M_{\la,\ka}$ in the derived category of $(\ghat_{\ka},\wt{B})$
Harish-Chandra modules, where $\wt{B}$ is the Iwahori subgroup of
$G[[t]]$, the preimage of a Borel subgroup $B \subset G$ under the
evaluation homomorphism $G[[t]] \to G$. This algebra is realized as
the cohomology
$$
H^\bullet(\g\ppart,\h,\on{End}_\C \M_{\la,\ka}) \simeq
H^\bullet(\tb,\h,\M_{\la,\ka} \otimes \C_{-\la}),
$$
where $\tb$ is the Lie algebra of $\wt{B}$ and $\h$ is the Cartan
subalgebra of the constant subalgebra $\g$. It follows from the
results of \cite{FF:gd,F:wak} that $\on{End}_{\ghat_{\ka_c}}
\M_{\la,\ka_c}$ is the algebra of functions on the space $\on{Op}_{^L
\g}(D)_\la$ of $^L \g$--opers on $D$ with regular singularities and
residue $-\la-\rho$. We show, in the same way as in the case of the
vacuum modules, that the full algebra of self-extensions of
$\M_{\la,\ka_c}$ is isomorphic to the algebra of differential forms on
$\on{Op}_{^L \g}(D)_\la$. We also show that for $\ka \neq \ka_c$ this
algebra if just isomorphic to $\C$.

The above statements are closely related to and were motivated by the
result of B. Feigin announced in \cite{Fe} that the cohomology of
$\ghat_{\ka_c}$ with coefficients in the completion of the enveloping
algebra of $\ghat_{\ka_c}$ (with respect to the adjoint action) is
isomorphic to the algebra of differential forms on the space of $^L
\g$--opers on the punctured disc.

The paper is organized as follows. In \secref{qc case} we consider the
graded version of the vacuum module (we call it ``classical'') and
recall the result of \cite{FGT}. In \secref{vac module} we define the
vacuum module $\V_{\ka}$ over the affine Kac-Moody algebra
$\ghat_{\ka}$ and describe the algebra of its self-extensions in a
suitable category of Harish-Chandra modules as the relative Lie
algebra cohomology $H^\bullet(\g\ppart,\g,\on{End}_\C \V_{\ka})$. In
\secref{DGVA} we realize it as the cohomology of the Chevalley complex
of $\g[[t]]$ relative to $\g$ with coefficients in $\V_{\ka}$. We show
that the algebra structure on this cohomology is induced by a DG
vertex superalgebra structure on the complex. Next, we show that the
latter is actually skew-commutative, and hence gives rise to the
structure of skew-commutative associative algebra on
$H^\bullet(\g[[t]],\V_{\ka_c})$. Our first proof of this fact was
based on an explicit computation of this cohomology. But, as was
subsequently pointed out to us by D. Gaitsgory, this is in fact a
corollary of a general property of the Chevalley complex of an
arbitrary vertex Lie superalgebra (see \propref{homotopy} and
\remref{generalization}).

In \secref{coh vac} we compute the cohomology of the vacuum module
$\V_{\ka_c}$ and show that it is isomorphic to the algebra of
differential forms on the space $\on{Op}_{^L \g}(D)$. We then show
that the cohomology of $\V_{\ka}$ where $\ka \neq \ka_c$ is isomorphic
to $\C$. In \secref{class verma} we consider the graded version
$\M^{\on{cl}}$ of the Verma module and compute the cohomology
$H^\bullet(\tb,\h,\M^{\on{cl}})$ using the results of \cite{FGT}.
Finally, in \secref{q verma} we use the computation in the graded case
and a description of the algebra of endomorphisms of the Verma modules
$\M_{\la,\ka}$ to find the full algebra of self-extensions of
$\M_{\la,\ka}$ in the appropriate category of Harish-Chandra
modules.

\medskip

\noindent{\bf Acknowledgments.} We wish to thank B. Feigin and
D. Gaitsgory for useful discussions. We are also grateful to
D. Gaitsgory for his careful reading of the manuscript and his
comments and suggestions.

The research of E.F. was supported by grants from DARPA and NSF, and
the research of C.T. was supported by grants from NSF and EPSRC.

\section{Cohomology of the classical vacuum module}    \label{qc case}

Let $\g$ be a simple finite-dimensional Lie algebra and $G$ the
connected simply-connected algebraic group corresponding to
$\g$. Introduce the ``classical'' vacuum module over the Lie algebra
$\g[[t]]$,
$$
\V^{\on{cl}} = \on{Sym}(\g\ppart/\g[[t]]) \simeq \on{Fun}(\g^*[[t]]
dt),
$$
where we use the residue pairing. We consider the relative cohomology
$H^\bullet(\g[[t]],\g,\V^{\on{cl}})$ of $\g[[t]]$ modulo the constant
Lie subalgebra $\g$ with coefficients in $\V^{\on{cl}}$. This
cohomology may be computed by the standard Chevalley complex of Lie
algebra cohomology
$$
C^\bullet(\g[[t]],\g,\V^{\on{cl}}) = \left( \V^{\on{cl}} \otimes
\bigwedge{}^\bullet(\g[[t]]/\g)^* \right)^{\g} = \left(
\on{Fun}(\g^*[[t]] dt) \otimes \bigwedge{}^\bullet(\g[[t]]/\g)^*
\right)^{\g}.
$$
\begin{remark}\label{topology}
Here and below we view $\g[[t]]$ and similar vector spaces as complete
topological spaces, and by the dual of such a vector space we mean the
vector space of continuous linear functionals on it (the topological
dual). As the result, $C^\bullet(\g[[t]],\g,\V^{\on{cl}})$ is a vector
space topologized as the direct limit of its finite-dimensional
subspaces.\qed
\end{remark}

The commutative algebra structure on $\V^{\on{cl}}$ gives rise to a
skew-commutative algebra structure on the complex
$C^\bullet(\g[[t]],\g,\V^{\on{cl}})$. The Chevalley differential is a
derivation of this algebra structure. Hence the cohomology
$H^\bullet(\g[[t]],\g,\V^{\on{cl}})$ also acquires the structure of a
skew-commutative algebra. We recall the description of this algebra
obtained in \cite{FGT}.

Set $$C_{\g^*} = \g^*/G := \on{Spec} (\on{Fun} \g^*)^G$$ and define the
local Hitchin space $C_{\g^*,\om}$ as
$$
C_{\g^*,\om} = \Gamma(D,\om \underset{\C^\times}\times C_{\g^*}),
$$
where $\om$ is the canonical line bundle on the disc $D = \on{Spec}
\C[[t]]$. Let us choose generators $P_i, i=1,\ldots,\ell$, of
$(\on{Fun} \g^*)^G$ of degrees $d_i+1$, where the $d_i$'s are the
exponents of $\g$, so that $(\on{Fun} \g^*)^G =
\C[P_i]_{i=1,\ldots,\ell}$. Then we obtain an identification
$$
\on{Fun}(C_{\g^*,\om}) = \C[P_{i,n}]_{i=1,\ldots,\ell;n \geq 0},
$$
where the $P_{i,n}$'s are the functions on $\g^*[[t]] dt$ defined by
the formula
\begin{equation}    \label{Pin}
P_{i,n}(f(t) dt) = \on{the} \; t^n\on{-coefficient} \; \on{of} \;
P_i(f(t)).
\end{equation}
Note that this definition depends on the choice of $t$, but we obtain
a coordinate-independent isomorphism
$$
C_{\g^*,\om} \simeq \bigoplus_{i=1}^\ell
\Gamma(D,\om^{\otimes(d_i+1)}) = \C[[t]] (dt)^{\otimes(d_i+1)}
$$
(it depends only on the choice of the generators $P_i$ in $(\on{Fun}
\g^*)^G$).

It is clear that the functions $P_{i,n}$ on $\g^*[[t]] dt$ are
$\g[[t]]$--invariant. Hence we obtain a homomorphism
$$
\on{Fun}(C_{\g^*,\om}) \to H^0(\g[[t]],\g,\V^{\on{cl}}).
$$

Next, following \cite{FGT}, we construct a map
$$
\varphi^{\on{cl}}_{\ka_0}: \on{Fun}(C_{\g^*,\om}) \to
H^1(\g[[t]],\g,\V^{\on{cl}}),
$$
for any non-zero invariant inner product $\kappa_0$ on $\g$ (recall
that such $\kappa_0$ is unique up to a scalar). By abusing notation, we
will also denote by $\kappa_0$ the corresponding map $\g^* \to \g$. We
need to associate to each $P \in \on{Fun}(C_{\g^*,\om})$ a cocycle
$\varphi^{\on{cl}}_{\ka_0}(P)$ in
$$
C^1(\g[[t]],\g,\V^{\on{cl}}) = \on{Hom}(\g[[t]]/\g,\on{Fun}(\g^*[[t]]
dt))^{\g}.
$$
We set
\begin{equation}    \label{varphi}
(\varphi^{\on{cl}}_{\kappa_0}(P))(x(t)) = \pa_{\kappa_0(dx(t))} \cdot
  P = \langle \kappa_0(dx(t)),dP\rangle,
\end{equation}
where $\kappa_0(dx(t)) \in \g^*[[t]] dt$ is considered as the constant
tangent vector field to $\g^*[[t]] dt$, and $\pa_{\kappa_0(dx(t))}$ is
the corresponding directional derivative. It is clear from this
formula that $\varphi^{\on{cl}}_{\ka_0}$ factors through the de Rham
differential $d: \on{Fun}(C_{\g^*,\om}) \to \Omega^1(C_{\g^*,\om})$.

In particular, applying the map $\varphi^{\on{cl}}_{\ka_0}$ to the
space $C_{\g^*,\om}^*$ of generators of $\on{Fun}(C_{\g^*,\om})$, which
is the dual space to $C_{\g^*,\om}$, we obtain a map
$$
C_{\g^*,\om}^* \to H^1(\g[[t]],\g,\V^{\on{cl}}).
$$
If we multiply $\ka_0$ by $\lambda$, then this map will get divided
by $\la$.

The following theorem is a version of Theorem B of \cite{FGT}.

\begin{thm}    \label{classical}
There is an isomorphism of graded algebras
$$H^\bullet(\g[[t]],\g,\V^{\on{cl}}) \simeq
\Omega^\bullet(C_{\g^*,\om}) = \on{Fun}(C_{\g^*,\om}) \otimes
\bigwedge{}^\bullet(C_{\g^*,\om}^*).$$ The right hand side is a free
skew-commutative algebra with even generators $P_{i,n} \in
H^0(\g[[t]],\g,\V^{\on{cl}})$ given by formula \eqref{Pin} and odd
generators $\varphi^{\on{cl}}_{\kappa_0}(P_{i,n}) \in$
$H^1(\g[[t]],\g,\V^{\on{cl}})$ given by formula \eqref{varphi}.
\end{thm}

We will now explain how to obtain the odd generators from the even
ones by using a deformation of the $\g[[t]]$--module $\V^{\on{cl}}$.

For $h \in \C$, let $\on{Conn}_h$ be the space of
$h$--connections on the trivial $G$--bundle on the disc $D$. These are
operators of the form $h \pa_t + A(t)$, where $A(t) dt$ is a section
of $\g \otimes \om$ on $D$. The Lie algebra $\g[[t]]$ acts on this
space by infinitesimal gauge transformations:
\begin{equation}    \label{gauge}
x(t) \cdot A(t) = [x(t),A(t)] - h \pa_t x(t).
\end{equation}
In particular, $\on{Conn}_h$ becomes the space $\g[[t]] dt$ when
$h=0$. We set
$$
\V_h^{\on{cl}} = \on{Fun}(\on{Conn}_h).
$$
Recall that $\V^{\on{cl}} = \on{Fun}(\g^*[[t]] dt)$. Using the
invariant inner product $\kappa_0$ on $\g$, we identify $\g^*$ with
$\g$ and hence $\V^{\on{cl}}$ with $\on{Fun}(\g[[t]] dt)$. Then
$\V_h^{\on{cl}}$ becomes a one-parameter deformation of
$\V^{\on{cl}}$. Note that the spaces $\on{Conn}_h$ are isomorphic to
each other for all non-zero values of $h$, and so are the modules
$\V_h^{\on{cl}}$.

Explicitly, the action \eqref{gauge} translates into the following
action of $\g[[t]]$ on $\V_h^{\on{cl}}$. For $A \in \g, n \in \Z$,
denote $A \otimes t^n$ by $A_n$. Then the difference between the
actions of $A_n \in \g[[t]]$ on monomial elements $B_{1,m_1} \ldots
B_{k,m_k}, m_i < 0$, in $\V^{\on{cl}}_h$ and $\V^{\on{cl}}$ is equal
to
\begin{equation}    \label{compare}
h \sum_{i=1}^k n \delta_{n,-m_i} \kappa_0(A,B_i) B_{1,m_1} \ldots
\wh{B}_{i,m_i} \ldots B_{k,m_k}.
\end{equation}

When $h \neq 0$, any $h$--connection can be brought to the form
$h \pa_t$ by using gauge transformations from the first congruence
subgroup
\begin{equation}    \label{G1}
G^{(1)} = \{ g \in G[[t]] \, | \, g(0) = 1 \}.
\end{equation}
Therefore the Lie algebra $\g \otimes t\C[[t]]$ acts co-freely on
$\V_h^{\on{cl}}$ and we obtain the following result.

\begin{lem}
For $h \neq 0$ the cohomology $H^i(\g[[t]],\g,\V_h^{\on{cl}})$
vanishes for $i>0$, and $H^0(\g[[t]],\g,\V_h^{\on{cl}}) = \C$.
\end{lem}

The one-parameter family of $\g[[t]]$--modules $\V_h^{\on{cl}}$ gives
rise to a $\g[[t]] \otimes \C[h]$--module, free over $\C[h]$, which by
abuse of notation we also denote by $\V_{{h}}^{\on{cl}}$. Consider
the Chevalley complex $C^\bullet(\g[[t]],\g,\V^{\on{cl}}_{h})$ of
relative cohomology of $\g[[t]]$ modulo $\g$ with coefficients in this
module. Given a class $\omega$ in $H^i(\g[[t]],\g,\V^{\on{cl}})$, we
choose a cocycle $\wt{\omega}$ representing it and extend it in an
arbitrary way to an element $\wt{\omega}({h})$ of
$C^\bullet(\g[[t]],\g,\V^{\on{cl}}_{h})$. Applying the differential of
this complex to $\wt{\omega}({h})$, dividing by $h$ and considering
the result modulo ${h}$, we obtain an element of
$C^{i+1}(\g[[t]],\g,\V^{\on{cl}})$. It is clear that this element is a
cocycle and that the corresponding cohomology class in
$H^{i+1}(\g[[t]],\g,\V^{\on{cl}})$ is independent of all
choices. Thus, we obtain well-defined linear maps
$$
\varphi_{\ka_0}^{\on{cl},i}:
H^i(\g[[t]],\g,\V^{\on{cl}}) \to H^{i+1}(\g[[t]],\g,\V^{\on{cl}})
$$
for all $i\geq 0$. Writing down the definition of
$\varphi_{\ka_0}^{\on{cl},0}$ explicitly, using formula
\eqref{compare}, we obtain that $\varphi_{\ka_0}^{\on{cl},0}$ acts on
$H^0(\g[[t]],\g,\V^{\on{cl}})$ precisely as the map
$\varphi^{\on{cl}}_{\kappa_0}$ given by formula \eqref{varphi}.

The differential of the Chevalley complex
$C^\bullet(\g[[t]],\g,\V^{\on{cl}}_{h})$ is an odd derivation of its
algebra structure. Therefore we obtain that the operators
$\varphi_{\ka_0}^{\on{cl},i}$ combine into an odd derivation of the
algebra $H^\bullet(\g[[t]],\g,\V^{\on{cl}})$, which we will also
denote by $\varphi^{\on{cl}}_{\kappa_0}$. We have
$\varphi_{\ka_0}^{\on{cl},i+1} \circ \varphi_{\ka_0}^{\on{cl},i} = 0$
for all $i \geq 0$.

Recall from \thmref{classical} that we have an isomorphism (dependent 
on $\ka_0$)
$$
H^\bullet(\g[[t]],\g,\V^{\on{cl}}) \simeq \Omega^\bullet(C_{\g^*,\om})
$$
where $\Omega^\bullet(C_{\g^*,\om})$ is the algebra of differential
forms on $C_{\g^*,\om}$. Under this identification the derivation
$\varphi^{\on{cl}}_{\kappa_0}$ becomes nothing but the de Rham
differential.

\section{Generalities on extensions}    \label{vac module}

Fix an invariant inner product $\kappa$ on $\g$ (recall that it is
unique up to a scalar), and let $\ghat_\ka$ denote the one-dimensional
central extension of $\g \otimes \C\ppart$,
$$
0 \to \C K \to \ghat_\ka \to \g \otimes \C\ppart \to 0
$$
with the commutation relations
\begin{equation}    \label{KM rel}
[A \otimes f(t),B \otimes g(t)] = [A,B] \otimes f(t) g(t) - (\kappa(A,B)
\on{Res} f dg) K,
\end{equation}
where $K$ is a central element. The Lie algebra $\ghat_\ka$ is the
affine Kac-Moody algebra associated to $\ka$.

Introduce the vacuum module of level $\ka$:
$$
\V_\ka = \on{Ind}_{\g[[t]] \oplus \C K}^{\ghat_{\ka}} \C,
$$
where $\g[[t]]$ acts on the one-dimensional space $\C$ by $0$, and $K$
as the identity. We denote by $v$ the generating vector of this module.

We wish to compute the algebra $\on{Ext}^\bullet(\V_\ka,\V_\ka)$ of
self-extensions of $\V_\ka$ in suitable categories of Harish-Chandra
modules (see Propositions \ref{HC exts} and \ref{HC ext} below), where
the multiplication is the Yoneda product. In this section we will show
that this algebra is given by the cohomologies of a suitable relative
Chevalley complex. In \secref{DGVA} we will show that this algebra is
skew-commutative. Then we will compute this algebra in \secref{coh
vac}.

In what follows by a representation of the group $G[[t]]$ we will
understand a direct limit of its finite-dimensional algebraic
representations. In particular, any vector $v$ in such a representation is
invariant under the $N$th congruence subgroup of $G[[t]]$ for
sufficiently large $N$ and hence satisfies $t^N\g[[t]] v = 0$. Let
$\on{Rep} G[[t]]$ be the category of all such representations.

Consider the category $HC(\ghat_{\ka},G[[t]])$ of {\em Harish-Chandra
modules for the pair \linebreak $(\ghat_{\ka},G[[t]])$}, i.e.,
$\ghat_\ka$--modules on which the action of $\g[[t]]$ exponentiates to
an action of $G[[t]]$, and such that $K$ acts as the identity. This is
a full subcategory within all $\ghat_{\ka}$--modules, closed under
forming kernels and cokernels, and is therefore abelian.

Next, we consider the relative Chevalley complex of continuous,
$\g$--invariant linear maps from $\bigwedge{}^\bullet \g\ppart/\g$ to
$\on{End}\V_\ka$. Here $\on{End}\V_\ka$ is topologized as
$$
\underset{\longleftarrow}\lim \on{Hom}(V_\al,\V_\ka), \quad \on{where}
\quad \V_\ka = \underset{\longrightarrow}\lim V_\al, \qquad \dim V_\al
< \infty.
$$
This complex is naturally a DG algebra, so its cohomology
$H^\bullet(\g\ppart,\g,\on{End}\V_\ka)$ acquires the structure of an
associative (super)algebra.

\begin{prop} \label{HC exts}
There is a natural isomorphism
\[
\on{Ext}^\bullet_{HC(\ghat_{\ka},G[[t]])}(\V_\ka,\V_\ka) \simeq
H^\bullet\left(\g\ppart,\g,\on{End} \V_\ka\right). 
\]
Moreover, the Yoneda product  corresponds to the cup-product on the
Chevalley complex.
\end{prop}

\noindent As $HC(\ghat_{\ka},G[[t]])$ does \textit{not} have enough
projectives or injectives, this requires an argument.


\begin{proof}
The cohomologies of the Chevalley complex of continuous, $G$--invariant 
linear maps $M\otimes \bigwedge{}^\bullet \g\ppart/\g \to N$ define the 
bi-additive functors 
$$
(M,N) \mapsto \on{E}^k(M,N).
$$
The $\on{E}^\bullet$ are equipped with functorial ``connecting 
homomorphisms", converting short exact sequences, in either variable, 
into long exact sequences. This situation is summarized by saying that 
the $\on{E}^\bullet$ form an \textit{exact connected right sequence of 
functors}, in each variable separately \cite[III \S2]{stro}, \cite
[IV \S7]{lang}. Note that $\on{E}^0 = \on{Hom}$, which is left exact. 

The Yoneda $\on{Ext}$'s, on the other hand, are the \textit{universal} 
exact connected right sequences starting with $\on{Hom}$ \cite[III \S1]
{stro}; they therefore map functorially to the $\on{E}^\bullet$. To show 
equality of the two functors, it suffices to check universality of 
the $\on{E}^\bullet$. This will follow from Grothendieck's \textit
{erasability} criterion \cite[Theorem 3.4.3]{stro}, but first we wish 
to spell out our map from $\on{Ext}$ to $\on{E}$.

Consider first $k=1$: a representative $N\to X(\eta)\to M$ of a class
$\eta\in \on{Ext}^1$ gives, after a choice of $G$--invariant linear
splitting, a crossed homomorphism $\g\ppart\times M \to N$. The latter
is the same as a one-cocycle in the relative Chevalley
complex. Changing the splitting modifies the cocycle by a coboundary,
so we get a well-defined Chevalley class. Note that our cocycle will
be continuous, because any linear splitting of $X(\eta)$ is
so.\footnote{It is not difficult to show that our Harish-Chandra
category is closed under forming extensions, whence it follows that
all one-cocycles are in fact continuous.} It is a standard exercise to
check additivity and functoriality of the resulting map $\on{Ext}^1
\to \on{E}^1$, along with the match of the first connecting maps in
the long exact sequence.

A Yoneda $k$--extension $N\to X^1\to X^2 \to\dots \to X^k \to M$
factors into a product of $1$-extensions, by an epic-monic
factorisation (i.e., representing each morphism $X^p \to X^{p+1}$ as
the composition $X^p \to \on{Im} X^p \to X^{p+1}$).  We chose
$G$--linear splittings of the short exact sequences
$$0\to \on{Im}X^{p-1} \to X^p \to \on{Im}X^p\to 0$$ to obtain
Chevalley one-cocycles, and assign to the extension $(X^\bullet)$ the
cup-product of the corresponding classes, obtained from the previous
construction. Clearly, that is independent of the splittings, but we
must check that this assignment descends to equivalence classes. That
is, for a morphism (commutative diagram) of $k$--extensions
\[
\begin{array}{ccccc}
N &\to &(X^\bullet) &\to &M \\
\downarrow & & \downarrow\dots\downarrow & & \downarrow \\
N &\to &(Y^\bullet) & \to & M
\end{array}
\]
the associated Chevalley classes agree. To see that, note that we can
construct commuting splittings of the two extensions (by compatible
splitting of each step in their epic-monic factorisations). The
associated Chevalley $(M,N)$--cocycles then agree, by commutativity of
the diagram; so our assignment of Chevalley \emph{classes} descends to
Yoneda equivalence classes.

We now have a transformation of functors which preserves the products,
and must check it is an isomorphism. We say that the functors $\on{E}^
\bullet(M,N)$ are \textit{erasable in $M$} if, for any $k>0$ and
$\eta\in \on{E}^k(M,N)$ there exists an epimorphism $f: M'\to M$ so
that $f^*\eta=0$ in $\on{E}^k(M',N)$. The $\on{Ext}^\bullet$ always
satisfy this condition, essentially by construction (take $M' = X^k$
above). Assuming erasability of $\on{E}^\bullet$, let us check that
our transformation $\Phi: \on{Ext}^\bullet \to \on{E}^\bullet$ is an
isomorphism, by induction on $k$; we know it for $k=0$. Let then
$\eta\in\ker\Phi: \on{Ext}^k(M,N) \to \on{E}^k(M,N)$ and assume that
$f^*\eta=0$ as above. With $K=\ker f$, it follows that $\eta
=\delta\zeta$ for some $\zeta\in\on{Ext}^{k-1} (K,N)$. Then,
$\delta\Phi(\zeta) = \Phi(\eta) = 0$, so $\Phi(\zeta)$ comes from
$\on{E}^{k-1}(M',N)$. By inductive assumption, the latter is
$\on{Ext}^{k-1}(M',N)$, so $\zeta$ comes from there and then $\eta =
\delta\zeta =0$. This shows injectivity of the map; surjectivity is
left to the reader.

Let us then verify erasability in $M$. Let $M' = \on{Ind}_{\g[[t]]
\oplus \C K}^{\ghat_{\ka}} M$ and note that restriction of cocycles
gives a map of Chevalley complexes
\begin{equation}\label{chevquiso}
C^\bullet(\g\ppart,\g; M',N) \to C^\bullet(\g[[t]],\g; M,N), 
\end{equation}
which we claim to be a quasi-isomorphism. Granting this for now, we
identify the right-hand Chevalley cohomologies above with the
$\on{Ext}$ groups in the category $\on{Rep} G[[t]]$. Indeed, observe
that $\on{Rep} G[[t]]$ contains enough injectives: for any
representation $N$ of $G$, the space $\Gamma(N)$ of global sections of
the vector bundle associated to $N$ over $G[[t]]/G$ is injective. Any
representation $V$ of $G[[t]]$ may be embedded into a representation
of this type, namely, into $\Gamma(V)$. Clearly, the Chevalley
cohomology of $\on{Hom}(M,\Gamma(N))$ coincides with
$\on{Ext}_{G[[t]]}(M,\Gamma(N)) = \on{Hom}^G(M,N)$ by de Rham's
theorem on the contractible space $G[[t]]/G$. Thus, we deduce the
isomorphism
\begin{equation}\label{shapiro}
\on{E}^k(M',N) \simeq \on{Ext}^k_{G[[t]]}(M,N).
\end{equation}

The functor $\on{Ext}^k$ is left erasable in $\on{Rep} G[[t]]$. Now,
if $M$ is already a representation of $\ghat_\ka$, the inclusion of
vector spaces $M\to M'$ that we used above to construct a map
\eqref{chevquiso} has a splitting $M' \to M$ which is a homomorphism
of $\ghat_\ka$--modules. This splitting lifts any class $\eta
\in\on{E}^k(M,N)$ to $\eta'\in \on{E}^k(M',N)$. The isomorphic image
of $\eta'$, via \eqref{shapiro}, in $\on{Ext}^k_{G[[t]]}(M,N)$
vanishes after further lifting to $\on{Ext}^k_{G[[t]]}(M_1,N)$ via an
epimorphism $M_1\to M$, which is an erasure of $\eta$ in $\on{Rep}
G[[t]]$ (by the above, we know that such a map $M_1 \to M$ exists). We
then have a commutative diagram
\[
\begin{array}{ccc}
\on{E}^k(M',N) & \overset{\sim}\to & \on{Ext}^k_{G[[t]]}(M,N) \\
\downarrow & & \downarrow \\
\on{E}^k(M'_1,N) & \overset{\sim}\to & \on{Ext}^k_{G[[t]]}(M_1,N)
\end{array}
\]
where $M'_1 = \on{Ind}_{\g[[t]] \oplus \C K}^{\ghat_{\ka}} M$, which
shows that the image of $\eta'$ under the vertical map is zero. Hence
the composition $\on{E}^k(M,N) \to \on{E}^k(M',N) \to \on{E}^k(M',N)$
obtained via the composite epimorphism $M_1' \to M' \to M$ sends
$\eta$ to zero as well. The epimorphism $M_1' \to M$ is then the
desired erasure of $\eta$ in $HC(\ghat_{\ka},G[[t]])$.

We are left to check that \eqref{chevquiso} is a quasi-isomorphism. 
Letting 
$$F_{-1}\g\ppart = 0,\quad F_0\g\ppart = \g[[t]], \quad F_1\g\ppart =
\g\ppart$$ gives increasing filtrations on $\bigwedge \g\ppart/\g$ and
$M'$, and hence on their tensor product. (Thus, $F_0 M' = M$ and $F_1
M' = M\oplus (t^{-1}\g[t^{-1}])\otimes M$.) This leads to a
complementary, descending filtration on the Chevalley complex for
$(\g\ppart,\g)$: the maps that vanish on the $p$th part of the tensor
product filtration live in $p$th degree. The associated graded complex
is $C^\bullet(\g[[t]],\g; K\otimes M,N)$, where $K$ is the Koszul
complex on $\g\ppart/\g[[t]]$. Consequently, the restriction
\eqref{chevquiso} gives an isomorphism of cohomologies, if we use the
Gr of the left-hand Chevalley complex. Because the filtration is
Hausdorff and complete on our space of continous linear maps, cocycles
and coboundaries can be lifted order-by-order to the complex from its
Gr, leading to an isomorphism of cohomologies in the original map
\eqref{chevquiso}.\footnote{This is, of course, the ``collapse of the
spectral sequence" argument in \cite[Theorem~1.5.4]{Fu}.}
\end{proof}

Finally, consider the category $HC(\ghat_{\ka},G)$ of
\textit{continuous Harish-Chandra modules for the pair
$(\ghat_{\ka},G)$}: these are the $\g$--integrable
$\ghat_{\ka}$--modules (where $\g \subset \g_\ka$ is the constant
subalgebra) with continuous action for the direct limit topology (see
Remark \ref{topology}), and on which $K\equiv 1$. Continuity amounts
to asking that every vector is annihilated by some subalgebra
$t^N\g[[t]]$.

Note that the category $HC(\ghat_{\ka},G[[t]])$ is a full subcategory
of $HC(\ghat_{\ka},G)$.

\begin{prop}    \label{HC ext}
There are natural isomorphisms
$$
\on{Ext}^k_{HC(\ghat_{\ka},G)}(\V_\ka,\V_\ka) \simeq
\on{Ext}^k_{HC(\ghat_{\ka},G[[t]])}(\V_\ka, \V_\ka).
$$
\end{prop}

\begin{proof}
We shall see below that we have a natural isomorphism
\begin{equation}    \label{isom}
\on{Ext}^k_{HC(\g[[t]],G)}(M,N)\simeq \on{Ext}^k_{G[[t]]}(M,N)
\end{equation}
for any objects $M,N$ of $\on{Rep} G[[t]]$. Assuming this, we can
repeat the argument in the proof of \propref{HC exts}, with the
category $\on{Rep} G[[t]]$ replaced by $HC(\g[[t]],G)$, to conclude
that the Chevalley complex $C^\bullet(\g\ppart,\g; M,N)$ computes
the $\on{Ext}$ in $HC(\ghat_{\ka},G)$.

To verify the isomorphism \eqref{isom}, we use the right adjoint
functor $\on{CoInd}$ to the forgetful functor from $G[[t]]$--modules
to $HC(\g[[t]],G)$.  Namely, $\on{CoInd}(N)$ is the space of
horizontal sections of the bundle $G[[t]]\times_G N$ over $G[[t]]/G$,
with flat connection determined from the $\g[[t]]$--action. It is easy
to see that the right derived functors $R^q\on{CoInd}$ are the de Rham
cohomologies with coefficients in the same bundle. We also note that
modules of the form $\on{Hom}^G(U\g[[t]];N)$ are injective objects in
the category $HC(\ghat_{\ka},G[[t]])$. Indeed, for any object $M$ of
this category $\on{Hom}(M,\on{Hom}^G(U\g[[t]];N)$ is isomorphic to
$\on{Hom}^G(M,N)$. The functor $\on{CoInd}$ takes
$\on{Hom}^G(U\g[[t]];N)$ to $\Gamma(N)$. For any object $M$ of
$\on{Rep} G[[t]]$, we have therefore a Grothendieck spectral sequence
\[
\on{Ext}^p_{G[[t]]}(M,R^q\on{CoInd}(N))\Rightarrow 
\on{Ext}^{p+q}_{HC(\g[[t]],G)}(M,N)
\]
However, if $N$ is already in $\on{Rep} G[[t]]$, $R^q\on{CoInd}(N) =0$
for $q>0$: this is because the flat $N$--bundle over $G[[t]]/G$
becomes isomorphic to the trivial $N$--bundle by a shearing map
$(\gamma,n) \mapsto (\gamma,\gamma n)$, and contractibility of the
space implies the vanishing of higher cohomology. The spectral
sequence degenerates to the desired isomorphism of the $\on{Ext}$'s.
\end{proof}

\section{DG vertex algebra structure and skew-commutativity}
\label{DGVA}

We now wish to compute the algebra $\on{Ext}^\bullet_
{HC(\ghat_{\ka},G)}(\V_\ka,\V_\ka)$. By \propref{HC exts}, it is
isomorphic to the algebra $H^\bullet\left(\g\ppart,\g,\on{End}
\V_\ka\right)$. Since $\V_\ka$ is an induced module, Shapiro's lemma
\cite[Theorem~1.5.4]{Fu} (see also the end of our proof of \propref{HC
exts}) implies the following statement.

\begin{lem}    \label{isom with exts}
We have a canonical isomorphism
$$
H^\bullet(\g[[t]],\g,\V_\ka) \simeq
H^\bullet(\g\ppart,\g,\on{End} \V_\ka).
$$
In particular, $H^\bullet(\g[[t]],\g,\V_\ka)$ acquires an algebra
structure via this isomorphism.
\end{lem}

\begin{remark}
A similar argument shows that $H^\bullet(\g[[t]],\g,\V_\ka)$ is
isomorphic to the relative semi-infinite cohomology
$H^{\frac{\infty}{2}+\bullet} (\ghat_{2\ka_c},\g,\V_\ka \otimes
\V_{2\ka_c-\ka})$, which is well-defined because $\V_\ka \otimes
\V_{2\ka_c-\ka}$ is a module of twice the critical level, as defined
in \secref{coh vac}.\qed
\end{remark}

Thus, we need to compute the cohomology $H^\bullet(\g[[t]],\g,\V_\ka)$. It
is realized as the cohomology of the relative  Chevalley complex
$$
C^\bullet(\g[[t]],\g,\V_\ka) = (\V_\ka \otimes
\bigwedge{}^\bullet (\g[[t]]/\g)^*)^\g.
$$
However, it is not immediately clear what algebra structure is induced
on the cohomology of this complex via the isomorphism of \lemref{isom
with exts}. We now explain how to define this algebra structure
directly on $H^\bullet(\g[[t]],\g,\V_\ka)$.

Namely, we will define a DG vertex superalgebra structure on the
relative Chevalley complex $C^\bullet(\g[[t]],\g,\V_\ka)$. It will induce
a vertex superalgebra structure on its cohomology
$H^\bullet(\g[[t]],\g,\V_\ka)$. We will show that the vertex superalgebra
$H^\bullet(\g[[t]],\g,\V_\ka)$ is skew-commutative, and hence we will
obtain the structure of an ordinary skew-commutative algebra on
$H^\bullet(\g[[t]],\g,\V_\ka)$. Finally, we will show that the latter
coincides with the algebra structure obtained from \lemref{isom with
exts}. Hence, by \propref{HC exts}, the algebra
$H^\bullet(\g[[t]],\g,\V_\ka)$ is isomorphic to the algebra
$\on{Ext}^\bullet_ {HC(\ghat_{\ka},G)}(\V_\ka,\V_\ka)$.

We start by recalling that the $\ghat_{\ka}$--module $\V_\ka$ is a
vertex algebra (see \cite{FB}, \S~2.4). Let $\{ J^a \}$ be a basis of
$\g$. Denote the element $J^a \otimes t^n \in \ghat_{\ka}$ by
$J^a_n$. Then $\V_\ka$ is freely generated by the vertex operators
$$
Y(J^a_{-1} v,z) = J^a(z) = \sum_{n \in \Z} J^a_n z^{-n-1}
$$
in the sense of the Reconstruction Theorem \cite{FB}, \S~4.4.

Next, recall that any skew-commutative associative superalgebra $V$
with a unit and a superderivation $T$ carries a canonical vertex
superalgebra structure with the vertex operators defined by the
formula \cite[\S~1.4]{FB}
$$
Y(A,z) = \on{mult}(e^{zT} A) = \sum_{n\geq 0} \frac{1}{n!} z^n
\on{mult}(T^n A).
$$
Conversely, given a skew-commutative vertex superalgebra $V$, with the 
vertex operation
\begin{equation}    \label{YY}
Y: V \to \on{End} V[[z^{\pm 1}]], \qquad A \mapsto Y(A,z) = \sum_{n
  \in \Z} A_{(n)} z^{-n-1},
\end{equation}
we recover the skew-commutative algebra product on it by the formula
$A,B \mapsto A_{(-1)} B$.

Consider the skew-commutative algebra
$\bigwedge{}^\bullet(\g[[t]]^*)$. It has generators $\psi^*_{a,n} \in
\g[[t]]^*$, where $a = 1,\ldots,\dim \g, n \leq 0$, defined by the
formula $\psi^*_{a,n}(J^b_m) = \delta_{a,b} \delta_{n,-m}$. Define a
superderivation $T$ on $\bigwedge{}^\bullet(\g[[t]]^*)$ by the formula
$T \cdot \psi^*_{a,n} = -(n-1) \psi^*_{a,n-1}$. Then
$\bigwedge{}^\bullet(\g[[t]]^*)$ acquires the structure of a vertex
superalgebra. It is freely generated by the vertex operators
$$
Y(\psi^*_{a,n},z) = \psi^*_a(z) = \sum_{n\leq 0} \psi^*_{a,n} z^{-n},
$$
in the sense of the Reconstruction Theorem \cite{FB}, \S~4.4.

Let us observe that $\bigwedge{}^\bullet(\g[[t]]^*)$ is a module over
the Clifford algebra associated to the vector space $\g[[t]] \oplus
\g[[t]]^*$ with the non-degenerate symmetric bilinear form induced by
the residue pairing. This algebra has generators $\psi_{a,n},
\psi^*_{a,m}, a = 1,\ldots,\dim \g; n \geq 0, m \leq 0$, satisfying
the anti-commutation relations
$$
[\psi_{a,n},\psi^*_{b,m}]_+ = \delta_{a,b} \delta_{n,-m}, \qquad
[\psi_{a,n},\psi_{b,m}]_+ = [\psi^*_{a,n},\psi^*_{b,m}]_+ = 0.
$$
The operators $\psi^*_{a,m}$ act on $\bigwedge{}^\bullet(\g[[t]]^*)$
by multiplication and the operators $\psi_{a,n}$ act by
contraction. Introduce the generating functions
$$
\psi_a(z) = \sum_{n \geq 0} \psi_{a,n} z^{-n-1}.
$$

We consider the tensor product vertex superalgebra structure on the
Chevalley complex
$$
C^\bullet(\g[[t]],\V_\ka) = \V_\ka \otimes \bigwedge{}^\bullet(\g[[t]]^*).
$$
For $A \in C^\bullet(\g[[t]],\V_\ka)$ we denote by $p(A)$ its parity. The vertex
operation
\begin{align*}
Y: C^\bullet(\g[[t]],\V_\ka) & \to \on{End} C^\bullet(\g[[t]],\V_\ka)[[z^{\pm
      1}]],
\end{align*}
may be determined by linearity by the following explicit formula
\begin{multline}    \label{Y}
Y(J^{a_1}_{n_1} \ldots J^{a_k}_{n_k} v \otimes \psi^*_{b_1,m_1} \ldots
\psi^*_{b_l,m_l},z) = \prod_{i=1}^k \frac{1}{(-n_i-1)!}  \prod_{j=1}^l
\frac{1}{(-m_j)!} \times \\ {\bf :} \pa_z^{-n_1-1} J^{a_1}(z) \ldots
\pa_z^{-n_k-1} J^{a_k}(z) {\bf :} \, \pa_z^{-m_1} \psi^*_{b_1}(z)
\ldots \pa_z^{-m_l} \psi^*_{b_l}(z),
\end{multline}
where the colons denote normal ordering.

The differential $d$ of Lie algebra cohomology is given by the formula
$$
d = \on{Res}_{z=0} \left( \sum_a J^a(z) \psi^*_a(z) - \frac{1}{2}
\sum_{a,b,c} \mu^{ab}_c \psi^*_a(z) \psi^*_b(z) \psi_c(z) \right) dz,
$$
where $(\mu^{ab}_c)$ are the structure constants of $\g$:
$$
[J^a,J^b] = \sum_c \mu^{ab}_c J^c.
$$
One checks easily that $d$ is a superderivation of the vertex algebra
$C^\bullet(\g[[t]],\V_\ka)$, i.e., we have
\begin{equation}    \label{DG}
Y(dA,z) = [d,Y(A,z)]_{\pm},
\end{equation}
where the sign of the commutator depends on whether the parity of $A
\in C^\bullet(\g[[t]],\V_\ka)$ is even or odd. Therefore the differential
$d$ gives $C^\bullet(\g[[t]],\V_\ka)$ the structure of a DG vertex
superalgebra. Thus, its cohomology is a graded vertex superalgebra.

Next, observe that the relative Chevalley complex
$C^\bullet(\g[[t]],\g,\V_\ka)$ is equal to the intersection of the kernels
of the operators $\psi_{i,0} = \on{Res}_{z=0} \psi_i(z) dz$ and
$\wh{J}^a_0 = \on{Res}_{z=0} \wh{J}^a(z) dz$, $a=1,\ldots,\dim \g$,
where
$$
\wh{J}^a(z) = J^a(z) - \sum_b \mu^{ab}_c \psi_b^*(z) \psi_c(z).
$$
These operators are superderivations of the vertex algebra structure
on $C^\bullet(\g[[t]],\V_\ka)$, i.e., they satisfy relations similar to
\eqref{DG}. Therefore the interesection of the kernels of these
operators is a vertex subalgebra of $C^\bullet(\g[[t]],\V_\ka)$. Hence we
obtain that $C^\bullet(\g[[t]],\g,\V_\ka)$ also carries the structure of a
DG vertex superalgebra, and so its cohomology
$H^\bullet(\g[[t]],\g,\V_\ka)$ is a graded vertex superalgebra.

The following result (and its generalization described in
\remref{generalization}) has been suggested to us by D. Gaitsgory.

\begin{prop}    \label{homotopy}
The vertex superalgebras $C^\bullet(\g[[t]],\V_\ka)$ and
$C^\bullet(\g[[t]],\g,\V_\ka)$ are homotopy skew-commutative. Thus, their
cohomology vertex superalgebras $H^\bullet(\g[[t]],\V_\ka)$ and
$H^\bullet(\g[[t]],\g,\V_\ka)$ are skew-commutative.
\end{prop}

\begin{proof}
Suppose that $V$ is a vertex superalgebra with the vertex operation
$Y$ (see formula \eqref{YY}).  Restricting it to the negative powers of
$z$ we obtain an operation $Y_-: V \to \on{End} V \otimes z^{-1}
\C[[z^{-1}]]$, which gives $V$ the structure of a vertex Lie
superalgebra (see \cite{FB}, \S~16.1). Recall that a vertex
superalgebra is called skew-commutative if $Y_- \equiv 0$.

We will introduce the following notation for $A, B \in V$,
\begin{align*}
Y_{(m)}(A,B) &= A_{(m)} B, \qquad m \geq 0, \\
Y_{(m)}(A,B) &= 0, \qquad m<0.
\end{align*}
Then the vertex Lie superalgebra operation $Y_-$ on $V$ is encoded by
the linear maps $Y_{(m)}: V \otimes V \to V$. Thus, $V$ is
skew-commutative if and only if $Y_{(m)} = 0$ for all $m \geq 0$.

We recall that the operations $Y_{(m)}$ satisfy the following
identities:
\begin{align}     \label{Ym1}
Y_{(m)}(A,B) &= (-1)^{p(A)p(B)} \sum_{n=0}^m \frac{1}{n!} (-1)^{m-n}
T^n Y_{(m-n)}(B,A), \\ \label{Ym2}
Y_{(m)}(TA,B) &= Y_{(m-1)}(A,B), \\ \label{Ym3}
Y_{(m)}(A,B_{(-1)} C) &= Y_{(m)}(A,B)_{(-1)} C + (-1)^{p(A) p(B)}
B_{(-1)} Y_{(m)}(A,C) \\ \notag &+ \sum_{j=0}^{m-1}
Y_{(m-j-1)}(Y_{(j)}(A,B),C)
\end{align}
(see \cite{FB}, \S~16.1).

Now consider the DG vertex superalgebra $C^\bullet =
C^\bullet(\g[[t]],\V_\ka)$ with the differential $d$. Formula \eqref{DG}
implies the following identity:
\begin{equation}    \label{homotopy identity}
d Y_{(m)}(A,B) = Y_{(m)}(dA,B) + (-1)^{p(A)} Y_{(m)}(A,dB).
\end{equation}

We will show that $C^\bullet$ is homotopy skew-commutative, i.e., we
will construct bilinear maps $$Z_{(m)}: C^\bullet \otimes C^\bullet
\to C^\bullet, \qquad m \geq 0,$$ of cohomological degree $-1$ such
that
\begin{equation}    \label{homotopy formula}
d Z_{(m)}(A,B) - Z_{(m)}(dA,B) + (-1)^{p(A)} Z_{(m)}(A,dB) =
Y_{(m)}(A,B).
\end{equation}
This will imply that the cohomology vertex superalgebra
$H^\bullet(\g[[t]],\V_\ka)$ is skew-commuta\-tive.

In order to construct it, we need to formulate some general properties
of biderivations of vertex Lie algebras. Let $V$ be a vertex Lie
algebra. We will call a collection of linear maps $Z_{(m)}: V \otimes
V \to V, m \geq 0$, a skew-symmetric biderivation if they satisfy the
following properties
\begin{align}     \label{Zm1}
Z_{(m)}(A,B) &= (-1)^{p(A)p(B)} \sum_{n=0}^m \frac{1}{n!} (-1)^{m-n}
T^n Z_{(m-n)}(B,A), \\ \label{Zm2}
Z_{(m)}(TA,B) &= Z_{(m-1)}(A,B), \\ \label{Zm3}
Z_{(m)}(A,B_{(-1)} C) &= Z_{(m)}(A,B)_{(-1)} C + (-1)^{p(A) p(B)}
B_{(-1)} Z_{(m)}(A,C) \\ \notag &+ \sum_{j=0}^{m-1}
Y_{(m-j-1)}(Z_{(j)}(A,B),C).
\end{align}
In particular, $(Y_{(m)})$ is a skew-symmetric biderivation.

The following lemma, which is proved by a straightforward calculation
using the identities \eqref{Ym1}--\eqref{Ym3},
\eqref{Zm1}--\eqref{Zm3} and \eqref{homotopy identity}, is a
generalization of the well-known properties of the ordinary
derivations.

\begin{lem} \hfill    \label{technical lemma}

\noindent {\em (1)} Let $V$ be a vertex superalgebra freely generated
by elements $A_i, i \in I$. Suppose that we are given
$Z_{(m)}(A_i,A_j) \in V$ for all $m \geq 0$ and $i,j \in I$,
satisfying the identities \eqref{Zm1}. Then this assignment may be
extended uniquely to a skew-symmetric biderivation of $V$.

\noindent {\em (2)} Let $V$ be a DG vertex superalgebra, freely
generated by $A_i, i \in I$, and $(Z_{(m)})$ a skew-symmetric
biderivation of $V$. Suppose that $(Z_{(m)})$ satisfies the identity
\eqref{homotopy formula} specialized to $A=A_i, B=A_j, i,j \in
I$. Then $(Z_{(m)})$ satisfies the identity \eqref{homotopy formula}
for all $A, B \in V$ and hence gives $V$ the structure of a homotopy
skew-commutative DG vertex superalgebra.

\end{lem}

Recall that our DG vertex superalgebra $C^\bullet$ is generated by
$J^a_{-1} v$ and $\psi^*_{a,0}, a =1,\ldots,\dim \g$. In light of
\lemref{technical lemma}, in order to construct a homotopy $(Z_{(m)})$
on $C^\bullet$, it is sufficient to define the $Z_{(m)}$'s on these
generators in such a way that they satisfy \eqref{homotopy formula}
and \eqref{Zm1}. We define them by the following formulas:
\begin{align*}
Z_{(m)}(J^a_{-1} v,\psi^*_{b,0}) &= - Z_{(m)}(\psi^*_{b,0},J^a_{-1} v)
= \frac{1}{2} \delta_{a,b} \delta_{m,0}, \\
Z_{(m)}(J^a_{-1} v,J^b_{-1} v) &= Z_{(m)}(\psi^*_{a,0},\psi^*_{b,0}) =
0.
\end{align*}

One checks by a direct calculation that these formulas satisfy the
identities \eqref{homotopy formula} and \eqref{Zm1}. Hence they can be
extended uniquely to a skew-symmetric biderivation of $C^\bullet$
which gives rise to the structure of a skew-commutative homotopy DG
superalgebra on it. This proves the statement of the proposition for
the Chevalley complex $C^\bullet = C^\bullet(\g[[t]],\V_\ka)$. In order to
prove the statement for the relative Chevalley complex
$C^\bullet(\g[[t]],\g,\V_\ka)$, which is a DG vertex subalgebra of
$C^\bullet$, it is sufficient to show that our biderivation
$(Z_{(m)})$ preserves this subcomplex. But it follows from the
construction that the maps $Z_{(m)}$ commute with the action of
$\psi^*_{a,0}$ and $\wh{J}^a_0$ on $C^\bullet$. Hence the homotopy
$(Z_{(m)})$ preserves $C^\bullet(\g[[t]],\g,\V_\ka)$. This completes the
proof.
\end{proof}

\begin{remark}    \label{generalization}
The vertex algebra $\V_\ka$ is the enveloping vertex algebra of the
vertex Lie algebra associated to the affine Kac-Moody algebra
$\ghat_{\ka}$. An analogue of the complex $C^\bullet(\g[[t]],\V_\ka)$
may be defined for the enveloping vertex algebra of any vertex Lie
algebra. This cohomology complex is then a DG vertex superalgebra. The
above proof carries over verbatim to this more general context, and we
obtain that this cohomology complex is always homotopy
skew-commutative.\qed
\end{remark}

Thus, we obtain that $H^\bullet(\g[[t]],\g,\V_\ka)$ is a skew-commutative
vertex superalgebra. Therefore the bilinear operation on
$H^\bullet(\g[[t]],\g,\V_\ka)$ given by the formula
\begin{equation}    \label{vertex product}
A, B \mapsto A_{(-1)} B
\end{equation}
defines the structure of a skew-commutative associative algebra on
it.

For a general vertex algebra this operation is non-commutative and
non-associative, so it does not give it the structure of an
algebra. But in our case, according to \propref{homotopy}, it gives
$H^\bullet(\g[[t]],\g,\V_\ka)$ the structure of a skew-commutative
associative algebra. In the next lemma we show that the above product
coincides with the product structure on the cohomology
$H^\bullet(\g[[t]],\g,\V_\ka)$ induced by its isomorphism with the algebra
of Ext's from \lemref{isom with exts}. Thus, we obtain that the latter
is skew-commutative, which is not obvious otherwise.

Let us observe that the Fourier coefficients of the series \eqref{Y}
may be viewed as elements of the Chevalley complex
$C^\bullet(\g\ppart,\g,\on{End} \V_\ka)$. Indeed, the $l$th group of this
complex is spanned by skew-symmetric continuous $l$--linear maps from
$\g\ppart$ to $\on{End} \V_\ka$. To a Fourier coefficient $A_{(n)} = \int
Y(A,z) z^n dz$ we assign the $l$--linear functional whose value on
$J^{b_1}_{p_1} \wedge \ldots \wedge J^{b_l}_{p_l}$ is the endomorphism
of $\V_\ka$ that equals the coefficient in front of $\psi^*_{b_1,p_1}
\ldots \psi^*_{b_l,p_l}$ in $A_{(n)}$.

\begin{lem}    \label{lem schapiro}
The isomorphism
\begin{equation}    \label{schapiro}
H^\bullet(\g[[t]],\g,\V_\ka) \simeq
H^\bullet(\g\ppart,\g,\on{End}_\C \V_\ka)
\end{equation}
of \lemref{isom with exts} may be realized as follows. Given a
cohomology class in $H^\bullet(\g[[t]],\g,\V_\ka)$, we represent it by a
cocycle $A$ in $C^\bullet(\g[[t]],\g,\V_\ka)$ and associate to it the
cohomology class of $A_{(-1)}$ which is a cocycle in
$C^\bullet(\g\ppart,\g,\on{End} \V_\ka)$. Then the natural product
structure on $H^\bullet(\g\ppart,\g,\on{End} \V_\ka)$ induces a product
structure on $H^\bullet(\g[[t]],\g,\V_\ka)$ which is given by formula
\eqref{vertex product}. In particular, the latter is automatically
associative.
\end{lem}

\begin{proof}
Formula \eqref{DG} implies that
$$
[d,A_{(-1)}]_\pm = (d A)_{(-1)}, \qquad A \in
C^\bullet(\g[[t]],\g,\V_\ka).
$$
Therefore if $A$ is a cocycle (resp., a coboundary) in
$C^\bullet(\g[[t]],\g,\V_\ka)$, then $A_{(-1)}$ is a cocycle (resp., a
coboundary) in $C^\bullet(\g\ppart,\g,\on{End} \V_\ka)$. Thus, we obtain a
well-defined map from $H^\bullet(\g[[t]],\g,\V_\ka)$ to
$H^\bullet(\g\ppart,\g,\on{End} \V_\ka)$. Let us show that this map
coincides with the isomorphism \eqref{schapiro} obtained via the
Shapiro lemma.

Recall the defining property of the isomorphism \eqref{schapiro}. If
$\omega$ is any cocycle in the Chevalley complex representing a
cohomology class $\ol{\omega}$ in $H^\bullet(\g\ppart,\g,\on{End} \V_\ka)$,
then the cohomology class in $H^\bullet(\g[[t]],\g,\V_\ka)$ corresponding
to $\ol{\omega}$ under the isomorphism \eqref{schapiro} is obtained by
restricting $\omega$ to $\bigwedge^\bullet(\g[[t]])$ and applying the
corresponding endomorphism of $\V_\ka$ to the vacuum vector $\V_\ka$. Recall
the vacuum axiom of the vertex algebra: in any vertex algebra $V$ with
the vacuum vector $\vac$ we have $A_{(-1)} \vac = A$ for any $A$. It
follows that if $\omega = A_{(-1)}$, then the result of the above
procedure will be $A$.

This proves that the isomorphism \eqref{schapiro} is indeed realized
by the assignment $A \mapsto A_{(-1)}$. But the product of $A_{(-1)}$
and $B_{(-1)}$ in $H^\bullet(\g\ppart,\g,\on{End}_\C \V_\ka)$, is just the
composition $A_{(-1)} \circ B_{(-1)}$. This induces the following
product structure on $H^\bullet(\g[[t]],\g,\V_\ka)$: $A,B \mapsto
(A_{(-1)} \circ B_{(-1)}) \vac = A_{(-1)} B$, i.e., the one given by
formula \eqref{vertex product}. This completes the proof.
\end{proof}

\section{Cohomology of the vacuum module of critical level}
\label{coh vac}

Now we specialize to the critical inner product $\ka_c = - \frac{1}{2}
\ka_{_K}$, where $\ka_{_K}$ denotes the Killing form on $\g$. Thus, by
definition, $$\ka_c(x,y) = - \frac{1}{2} \on{Tr} (\on{ad} x \on{ad}
y).$$

In order to simplify notation, we will denote $\V_{\ka_c}$ by $\V$.

In cohomological degree zero we have
$$
H^0(\g[[t]],\g,\V) = \V^{\g[[t]]} \simeq
\on{End}_{\ghat_{\ka_c}} \V.
$$
This algebra has been described in \cite{FF:gd,F:wak}. Let us recall
this result.

First we need to define $\g$--opers (see \cite{DS,BD,F:wak}). Choose a
Cartan decomposition $\g = \n_- \oplus \h \oplus \n$ and Chevalley
generators $f_i, i=1,\ldots,\ell$ of $\n_-$. Denote by $\bb$ the
direct sum $\h \oplus \n$ and by $N$ the Lie group of $\n$. Then by
definition a $\g$--oper on the disc $D = \on{Spec} \C[[t]]$ is an
equivalence class of first order operators
\begin{equation}    \label{another form of nabla}
\nabla = \pa_t + \sum_{i=1}^\ell f_i + {\mb v}(t), \qquad {\mb v}(t)
\in \bb[[t]],
\end{equation}
with respect to the gauge action of the group $N[[t]]$ given by the
formula
$$
g \cdot (\pa_t + A(t)) = \pa_t + g A(t) g^{-1} - \pa_t g \cdot g^{-1}.
$$
We denote the space of $\g$--opers on $D$ by $\on{Op}_\g(D)$.

Let us include the element $p_{-1} = \sum_{i=1}^\ell f_i$ into a
principal $\sw_2$--triple $\{ p_{-1},2\check{\rho},p_1 \}$, where
$\check\rho \in \h$ is the sum of the fundamental coweights. The space
of invariants in $\n$ of the adjoint action of $p_1$ has a basis $p_j,
j=1,\ldots,\ell$, of elements such that $[\check\rho,p_j] = d_j p_j$,
where $\{ d_1,\ldots,d_\ell \}$ is the set of exponents of $\g$.

It is known from \cite{DS} that the above action of $N[[t]]$ is free
and each equivalence class contains a unique operator of the form
$\pa_t + \sum_{i=1}^\ell f_i + {\mathbf v}(t)$, where
$${\mathbf v}(t) = \sum_{j=1}^\ell v_j(t) \cdot p_j, \qquad v_j(t) \in
\C[[t]].$$ The series $v_1(t)$ transforms as a projective connection,
while $v_j(t), j>1$, transforms as a $(d_j+1)$--differential (see,
e.g., \cite{F:wak}). Thus, we obtain an isomorphism
\begin{equation}    \label{repr}
\on{Op}_\g(D) \simeq  {\mc P}roj(D) \times \bigoplus_{j=2}^\el
\Gamma(D,\omega^{\otimes(d_j+1)})
\end{equation}
This isomorphism is coordinate-independent, but it depends in the
obvious way on the choice of the $f_i$'s and the $p_j$'s. The space
${\mc P}roj(D)$ is an affine space modeled on the space
$\Gamma(D,\omega^{\otimes 2})$ of quadratic differentials on $D$. It
follows from the identification \eqref{repr} that the algebra
$\on{Fun}(\on{Op}_{^L \g}(D))$ has a natural filtration, and the
associated graded algebra is canonically isomorphic to
$\on{Fun}(C_{\g^*,\om})$ (see \cite{BD}, \S\S 3.1.12--3.1.14, for more
details).

The module $\V$ has a PBW filtration, and the associated
graded $\on{gr} \V$ is isomorphic to $\V^{\on{cl}}$. This
induces a filtration on $\V^{\g[[t]]} =
\on{End}_{\ghat_{\ka_c}} \V$ and we obtain a natural
homomorphism of algebras
$$
\on{gr} \on{End}_{\ghat_{\ka_c}} \V \to (\V^{\on{cl}})^{\g[[t]]} =
\on{Fun}(C_{\g^*,\om}).
$$
Let $^L \g$ be the Lie algebra that is Langlands dual to $\g$, so that
its Cartan matrix is the transpose to that of $\g$.

\begin{thm}[\cite{FF:gd,F:wak}]    \label{center}
There is a canonical isomorphism of filtered algebras
$$
\on{End}_{\ghat_{\ka_c}} \V \simeq \on{Fun}(\on{Op}_{^L \g}(D))
$$
such that the diagram
$$
\CD
\on{gr} \on{End}_{\ghat_{\ka_c}} \V @>{\sim}>> \on{gr}
\on{Fun}(\on{Op}_{^L G}(D)) \\
@VVV    @VVV   \\   
\on{Fun}(C_{\g^*,\om}) @>{(-1)^{\on{deg}}}>> \on{Fun}(C_{\g^*,\om})
\endCD
$$
(where $(-1)^{\on{deg}}$ is the automorphism taking value $(-1)^n$ on
the subspace of elements of degree $n$) is commutative.
\end{thm}

Next, we construct maps $H^i(\g[[t]],\g,\V) \to
H^{i+1}(\g[[t]],\g,\V)$ similarly to the classical case, by
deforming the critical inner product. Consider the one-parameter
family of vacuum modules
$$
\V_\kappa = \on{Ind}_{\g[[t]] \oplus \C K}^{\ghat_{\ka}} \C,
$$
where $\kappa = \ka_c + h \ka_0$, with respect to the parameter
$h$. This family gives rise to a one-parameter family of Chevalley
complexes $C^\bullet(\g[[t]],\g,\V_\kappa)$.

\begin{prop}    \label{gen ka}
For $\ka \neq \ka_c$ the cohomology $H^i(\g[[t]],\g,\V_\ka)$ vanishes
for $i>0$, and $H^0(\g[[t]],\g,\V_\ka) = \C$.
\end{prop}

\begin{proof}
Consider the contragredient module $\V_\ka^\vee$. It follows from the
results of \cite{KK} that for generic $\kappa$ the $\ghat_\ka$--module
$\V_\ka$ is irreducible. Therefore the natural homomorphism of
$\ghat_\ka$--modules $\V_\ka \to \V_\ka^\vee$ is an isomorphism. It
follows that $\V_\ka$ is $\g \otimes t \C[[t]]$--cofree, and we obtain
the assertion of the proposition using Shapiro's lemma. A proof that
works for an arbitrary $\ka \neq \ka_c$ (when $\V_\ka$ may be
reducible) will be given after the proof of \thmref{vacuum}.
\end{proof}

Now, in the same way as in the classical case, we obtain linear maps
$$
\varphi^i_{\ka_0}: H^i(\g[[t]],\g,\V) \to
H^{i+1}(\g[[t]],\g,\V).
$$
In particular, for $i=0$ we obtain a map
$$
\varphi_{\ka_0} = \varphi^0_{\ka_0}: \on{Fun}(\on{Op}_{^L \g}(D)) \to
H^1(\g[[t]],\g,\V).
$$
Explicitly, this map looks as follows. Suppose that $P \in
\on{Fun}(\on{Op}_{^L \g}(D)) \simeq \V^{\g[[t]]}$ and let us compute
the corresponding one-cocycle in the Chevalley complex
$C^\bullet(\g[[t]],\g,\V)$, which is a linear map $\varphi_{\ka_0}(P):
\g \otimes t\C[[t]] \to \V$. It follows from the definition that
\begin{equation}    \label{varphi quantum}
x(t) \mapsto \left. \frac{1}{h} x(t) \cdot \wt{P} \right|_{h=0},
\end{equation}
where $\wt{P}$ is an arbitrary lifting of $P$ to $\V_{\kappa}$
considered as a free module over $\C[h]$.

\begin{lem}    \label{leibniz}
The maps $\varphi^i_{\ka_0}$ satisfy the Leibniz rule with respect to
the natural action of $\on{Fun}(\on{Op}_{^L \g}(D)) \simeq
\on{End}_{\ghat_{\ka_c}} \V$ on
$H^1(\g[[t]],\g,\V)$. In particular, $\varphi_{\ka_0}^0$
necessarily factors through the de Rham differential
$$
d: \on{Fun}(\on{Op}_{^L \g}(D)) \to \Omega^1(\on{Op}_{^L \g}(D)).
$$
We also have $\varphi_{\ka_0}^{i+1} \circ \varphi_{\ka_0}^i = 0$ for
all $i \geq 0$.
\end{lem}

\begin{proof}
After trivialising the terms $C^i(\g[[t]],\g,\V_\ka)$ of the Chevalley 
complex as bundles over $\on{Spec}\C[h]$, we can expand the differential 
$\delta_h$ in powers of $h$:
$$
\delta_h = \delta_0 + h \delta_1 + h^2 \delta_2 + \ldots.
$$
Note that $\delta_0$ is the differential on
$C^i(\g[[t]],\g,\V)$, and $\delta_1 = \varphi_{\ka_0}$ (here
we omit the upper index $i$ to simplify notation). From $\delta_h^2 =
0$, we find that $\delta_0^2 = 0$---which is clear anyway, as
$\delta_0$ is a differential; additionally, with $[\,,]_+$ denoting
the supercommutator,
$$
[\delta_0,\delta_1]_+ = 0\qquad\mbox{and}\qquad \delta_1^2 = 
-[\delta_0,\delta_2]_+.
$$
The first equality implies that $\delta_1 = \varphi_{\ka_0}$ maps
$\delta_0$--cocycles to $\delta_0$--cocycles and so indeed gives rise
to a well-defined map $H^i(\g[[t]],\g,\V) \to
H^{i+1}(\g[[t]],\g,\V)$, while the second equality means that
at the level of cohomologies we have $\delta_1^2 = 0$.

To see that the Leibniz rule holds, observe that $\delta_h$ preserves
the multiplicative structure on the Chevalley complex computing
$$H^\bullet(\g\ppart,\g,\on{End} \V_\ka) \simeq
H^\bullet(\g[[t]],\g,\V_\ka).$$ Thus we have
$$
\delta_h(A \underset{h}* B) = (\delta_h A) \underset{h}* B +
(-1)^{\deg A} A \underset{h}* (\delta_h B),
$$
where $\underset{h}*$ denotes the product on the complex depending on
the parameter $h$. Expanding $\delta_h, A$ and $B$ in powers of
$h$ and using the fact that $\delta_0$ preserves the multiplicative
structure on the Chevalley complex at $h=0$, we obtain that
$$
\delta_1(A_0 * B_0) = (\delta_1 A_0) * B + (-1)^{\deg A} A * (\delta_1
B).
$$
In particular, if we choose $A_0$ to be an element of
$H^0(\g\ppart,\g,\on{End} \V) = \on{Fun}(\on{Op}_{^L \g}(D))$, we
obtain the assertion of the lemma.
\end{proof}

Since $\on{Op}_{^L \g}(D)$ is an affine space modeled on the vector
space $C_{\g^*,\om}$, the map $\Omega^1(\on{Op}_{^L \g}(D)) \to
H^1(\g[[t]],\g,\V)$ gives rise to a map
\begin{equation}    \label{map}
C_{\g^*,\om}^* \to H^1(\g[[t]],\g,\V).
\end{equation}
If we multiply $\ka_0$ by $\lambda$, then the map \eqref{map} will get
divided by $\la$, as in the classical case.

Now we are ready to describe the cohomology algebra
$H^\bullet(\g[[t]],\g,\V)$. The following assertion means that this
algebra is a skew-commutative polynomial algebra with the even generators
being the generators of the polynomial algebra $\on{Fun}(\on{Op}_{^L
G}(D))$ and the odd generators being the images of some basis elements
of $C_{\g^*,\om}^*$ under the map \eqref{map} which turns out to be
injective.

\begin{thm}    \label{vacuum}
For each non-zero invariant inner product $\kappa_0$ on $\g$ we have
an isomorphism
\begin{equation}    \label{isom coh}
H^\bullet(\g[[t]],\g,\V) \simeq \Omega^\bullet(\on{Op}_{^L \g}(D)) =
\on{Fun}(\on{Op}_{^L \g}(D)) \otimes
\bigwedge{}^\bullet(C_{\g^*,\om}^*).
\end{equation}
where $\Omega^\bullet(\on{Op}_{^L \g}(D))$ is the algebra of
differential forms on $\on{Op}_{^L \g}(D)$, such that the map
$\varphi_{\ka_0}$ gets identified with the de Rham differential on
$\Omega^\bullet(\on{Op}_{^L \g}(D))$. If we rescale $\ka_0$ by
$\lambda$, then this isomorphism gets rescaled by $\la^{-i}$ on the
$i$th cohomology.
\end{thm}

\begin{proof}
Let us describe the strategy of the proof. We will compute the
cohomology $H^\bullet(\g[[t]],\g,\V)$ using the spectral
sequence induced by the PBW filtration on $\V$. The first term
of this spectral sequence is equal to the cohomology
$H^\bullet(\g[[t]],\g,\V^{\on{cl}})$. We will prove that all higher
differentials vanish. This will imply that, for the PBW filtration,
$$
\on{gr} H^\bullet(\g[[t]],\g,\V) =
H^\bullet(\g[[t]],\g,\V^{\on{cl}}),
$$
which is described by Theorem \ref{classical}. We will see that the
algebra structures agree and will verify the equality $\on{gr}\varphi_
{\ka_0} = \varphi_{\ka_0}^{\on{cl}}$. This will imply that
$H^\bullet(\g[[t]],\g,\V)$ is generated as an algebra in
degrees $0$ and $1$.  Theorem \ref{center} and Lemma \ref{leibniz}
will allow us to identify $\varphi_{\ka_0}$ with de Rham's
differential to give a linear map
$$
H^1(\g[[t]],\g,\V) \to \Omega^1(\on{Op}_{^L \g}(D))
$$
whose $\on{gr}$ is the isomorphism of Theorem \ref{classical}.
Symmetrized multiplication\footnote{we know from \propref{homotopy}
that $H^\bullet(\g[[t]],\g,\V)$ is commutative, but we
prefer not to use it here} of the generators will give us a linear map
in \eqref{isom coh}, whose associated graded is the isomorphism in
Theorem \ref{classical}. The theorem will then follow, since we know
from \propref{homotopy} and \lemref{lem schapiro} that the algebra
$H^\bullet(\g[[t]],\g,\V)$ is skew-commutative, and so
\eqref{isom coh} is an algebra isomorphism.

Now we proceed with the proof. First note that
$H^0(\g[[t]],\g,\V)$ is isomorphic to $\on{Fun}(\on{Op}_{^L
\g}(D))$ according to \thmref{center}. Since $\on{gr}
\on{Fun}(\on{Op}_{^L \g}(D)) = \on{Fun}(C_{\g^*,\om})$, we find from
\thmref{classical} that the entire zeroth cohomology part of the first
term of the spectral sequence survives. (Here and below we say that a
class in a particular term of a spectral sequence ``survives'' if it
is annihilated by all higher differentials of the spectral sequence
and does not lie in the image of any of the higher differentials.)

Now let us look at first cohomology part of the spectral
sequence. Recall that
$$
\on{Fun}(C_{\g^*,\om}) \simeq \C[P_{i,n}]_{i=1,\ldots,\ell;n\geq 0}.
$$
Let us choose elements $\wt{P}_{i,n}$ in the algebra
$\on{Fun}(\on{Op}_{^L \g}(D))$ whose symbols are equal to
$P_{i,n}$. Then we have
$$
\on{Fun}(\on{Op}_{^L \g}(D)) \simeq
\C[\wt{P}_{i,n}]_{i=1,\ldots,\ell;n\geq 0}.
$$
Applying the map $\varphi_{\ka_0}$ to the $\wt{P}_{i,n}$'s, we
obtain classes in the first cohomology group
$H^1(\g[[t]],\g,\V)$. To show that these classes are non-zero,
we compute how the symbol of $\varphi_{\ka_0}$ acts on the associate
graded cohomology and find that it coincides with the operator
$\varphi_{\ka_0}^{\on{cl}} = \varphi_{\ka_0}^{\on{cl},0}$.

Indeed, let us compute how the action of $A_n = A \otimes t^n \in
\g[[t]]$ on $\V$ changes when we deform the central
extension. We find that when we apply this element to a
lexicographically ordered monomial $B_{1,m_1} \ldots B_{k,m_k} v_{\ka}
\in \V_\ka$, the deformation is equal to
\begin{equation}    \label{differs}
h \sum_{i=1}^k n \delta_{n,-m_i} \kappa_0(A,B_i) B_{1,m_1} \ldots
\wh{B}_{i,m_i} \ldots B_{k,m_k} v_\ka
\end{equation}
plus the sum of monomials of order less than $k-1$. Thus, the
$h$--linear term of the deformation sends the $k$th term of the PBW
filtration on $\V$ to the $(k-1)$st term, and comparing with
formula \eqref{compare}, we find that the corresponding operator on
the associated graded coincides with the operator
$\varphi_{\ka_0}^{\on{cl}}$.

Therefore we find that the symbol of the cohomology class
$\varphi_{\ka_0}(\wt{P}_{i,n})$ is equal to
$\varphi_{\ka_0}^{\on{cl}}(P_{i,n})$, which is one of the odd
generators of $H^\bullet(\g[[t]],\g,\V^{\on{cl}})$ from
$H^1(\g[[t]],\g,\V^{\on{cl}})$. Thus, we obtain that the cocycles
representing all odd generators of
$H^\bullet(\g[[t]],\g,\V^{\on{cl}})$ may be lifted to cocycles in the
Chevalley complex $C^\bullet(\g[[t]],\g,\V)$. Since
$H^\bullet(\g[[t]],\g,\V)$ is computed by a spectral sequence whose
first term is $H^\bullet(\g[[t]],\g,\V^{\on{cl}})$ we find that these
cocycles can become trivial only if their symbols could be obtained
from some cohomology classes in $H^0(\g[[t]],\g,\V^{\on{cl}})$ under
the action of the higher differentials on our spectral sequence. But
we know that all classes in $H^0(\g[[t]],\g,\V^{\on{cl}})$ survive, so
under the action of the higher differentials they cannot kill any
classes in $H^1(\g[[t]],\g,\V^{\on{cl}})$. Hence all of these classes
survive, and we obtain an injective map $C_{\g^*,\om}^* \to
H^1(\g[[t]],\g,\V)$.

The natural action of $H^0(\g[[t]],\g,\V) =
\on{End}_{\ghat_{\ka_c}} \V$ on $H^i(\g[[t]],\g,\V)$
leads to a map
$$
\on{Fun}(\on{Op}_{^L \g}(D)) \otimes
C_{\g^*,\om}^* \to H^1(\g[[t]],\g,\V).
$$ 
In the same way as above, we find that it is injective. Considering 
symbols and comparing with $H^1(\g[[t]],\g,\V^{\on{cl}})$, we find
that this map is an isomorphism.

We now generate higher cohomology classes in
$H^i(\g[[t]],\g,\V), i>1$, by taking symmetrized products of
the cocycles representing the classes in $H^1(\g[[t]],\g,\V)$
that we have just constructed. We need to prove inductively that the
resulting cohomology classes are non-trivial and span
$H^\bullet(\g[[t]],\g,\V)$. In light of the classical result,
\thmref{classical}, it suffices to show that the multiplicative
structures on $H^\bullet(\g[[t]],\g,\V)$ and
$H^\bullet(\g[[t]],\g,\V^{\on{cl}})$ are compatible, in the sense that
the latter is the associated graded of the former.

We will give two proofs of this assertion. The first one is to use
formulas \eqref{vertex product} and \eqref{Y} defining the product
structure on $H^\bullet(\g[[t]],\g,\V)$ using the structure of
vertex algebra. Then we find immediately from the vacuum axiom of
vertex algebra that the symbol of the product of two cocycles $A, B
\in C^\bullet(\g[[t]],\g,\V)$, i.e., $\on{symb}(A_{(-1)} B)$
is equal to the product of their symbols, $\on{symb}(A) \on{symb}(B)$.

Now we give another proof that does not use the vertex superalgebra
structure on $C^\bullet(\g[[t]],\g,\V)$. Let us introduce a Lie
algebra filtration $(F_i)$ on $\g\ppart$ by setting
$$
F_{-1} = 0, \qquad F_0 = \g[[t]], \qquad F_1 = \g\ppart.
$$
Then the associated graded Lie algebra is
$$
\wt{\g} = \g[[t]] \oplus \g\ppart/\g[[t]],
$$
where the Lie algebra structure on the first summand is the usual one,
on the second summand it is commutative, and the commutator between
elements of the two summands is given by the natural action of
$\g[[t]]$ on $\g\ppart/\g[[t]]$. The induced filtration on
$\V$ coincides with the PBW filtration. Hence we have an
isomorphism
$$
\on{gr} \V \simeq \on{Ind}_{\g[[z]]}^{\wt\g} \C \simeq
\V^{\on{cl}}.
$$
Thus, using the Shapiro lemma, we obtain the isomorphism
$$
H^\bullet(\g[[t]],\g,\V^{\on{cl}}) \simeq
H^\bullet(\wt\g,\g,\on{End} \V^{\on{cl}}).
$$
which is compatible with the multiplicative structure. But the product
on the algebra $H^\bullet(\wt\g,\g,\on{End} \V^{\on{cl}})$ is clearly
the associated graded of the product on $H^\bullet(\g\ppart,\g,\on{End}
\V)$, and hence we obtain the desired assertion.

Finally, it follows from \propref{homotopy} and \lemref{lem schapiro}
that the algebra $H^\bullet(\g[[t]],\g,\V)$ is
skew-commutative. This completes the proof.
\end{proof}

\begin{remark}    \label{alternaive}
The skew-commutativity of $H^\bullet(\g[[t]],\g,\V)$ may be
shown without using the vertex superalgebra structure on
$C^\bullet(\g[[t]],\g,\V)$. By our construction,
$H^\bullet(\g[[t]],\g,\V)$ is generated by classes in $H^0$
and $H^1$, and we already know that all elements
$H^0(\g[[t]],\g,\V)$ come from central elements. Thus, $H^0$
is central in the entire cohomology, and it remains to show that the
generators in $H^1$ anti-commute.  These generators have the form
$\varphi_{\ka_0}(A)$ where $A \in H^0(\g[[t]],\g,\V)$. Their
anti-commutativity follows from \lemref{leibniz}. Indeed, we have
(denoting $\varphi_{\ka_0}$ by $\delta_1$ as in the proof of
\lemref{leibniz})
$$
[\delta_1(A),\delta_1(B)]_+ = \delta_1([A,\delta_1(B)]) - A
\delta_1^2(B) + \delta_1^2(B) A,
$$
using the Leibniz rule. The first term on the right hand side vanishes
because $A$ is central in $H^\bullet(\g[[t]],\g,\V)$, and the
remaining two terms vanish because $\delta_1^2 = 0$ on the
cohomology. Hence $[\delta_1(A),\delta_1(B)]_+ = 0$.\qed
\end{remark}

\medskip

\begin{proof}[Proof of \propref{gen ka}.] Consider the spectral
sequence induced by the PBW filtration on $\V_\ka$ and converging to
$H^\bullet(\g[[t]],\g,\V_\ka)$. It has the same first term
$H^\bullet(\g[[t]],\g,\V^{\on{cl}})$ as the one for $\V$, considered
in the proof of \thmref{vacuum}. The computation of $\varphi_{\ka_0}$
in the above proof shows that the first differential of this spectral
sequence for $\V_\ka$ differs from the one for $\V$ by the map
$h\varphi_{\ka_0}^{\on{cl}}$, where $\ka=\ka_0+h\ka_0$ (see formula
\eqref{differs}). But the differential for $\V$ was equal to zero,
and we know from \thmref{classical} that $\varphi_{\ka_0}^{\on{cl}}$
is the de Rham differential on $$H^\bullet(\g[[t]],\g,\V^{\on{cl}})
\simeq \Omega^\bullet(C_{\g^*,\om}).$$ Hence its cohomology is $\C$,
in degree $0$.
\end{proof}

\medskip

Next, we consider the absolute cohomology. 

\begin{prop}
The cohomology $H^\bullet(\g[[t]],\V)$ is canonically isomorphic to
the tensor product of $H^\bullet(\g[[t]],\g,\V)$ described in
\thmref{vacuum} and $H^\bullet(\g,\C)$ which is the exterior algebra
with generators in degrees $2d_i+1, i=1,\ldots,\ell$.
\end{prop}

\begin{proof}
The argument is standard, applying to any semi-direct product Lie
algebra $\g\ltimes\mathfrak{L}$ ($\mathfrak{L} = \g \otimes t\C[[t]]$
here), provided that $\g$ is reductive, $\mathfrak{L}$ is a direct
product, and $\V$ a direct sum of finite-dimensional irreducible
$\g$-modules.

Consider the Serre-Hochshild spectral sequence corresponding to the 
Lie subalgebra $\g \subset \g[[t]]$ (see, e.g., \cite{Fu}). In the
first term we have
$$
E_1^{p,q} = H^q\left(\g,\V \otimes \bigwedge{}^p(\g[[t]]/\g)^*\right) 
	= H^q(\g,\C) \otimes C^p(\g[[t]],\g,\V),
$$
where
$$
C^p(\g[[t]],\g,\V) = \left( \V \otimes \bigwedge{}^p(\g[[t]]/\g)^*
\right)^{\g}
$$
is the $p$th group of the relative Chevalley complex computing
$H^\bullet(\g[[t]],\g,\V)$. Therefore in the second term we have
\begin{equation}    \label{E2}
E_2^{p,q} = H^q(\g,\C) \otimes H^p(\g[[t]],\g,\V).
\end{equation}
We can therefore represent classes in $E_2^{p,q}$ as tensor products
$\omega_1 \otimes \omega_2$ of a cocycle $\omega_1$ in $C^q(\g,\C)$
representing a class in $H^q(\g,\C)$ and a cocycle $\omega_2$ in $(\V
\otimes \bigwedge^p(\g \otimes t\C[[t]])^*)^{\g}$ representing a
class in $H^p(\g[[t]],\g,\V)$. Applying the differential to this
class, we find that it is identically equal to zero because $\omega_2$
is $\g$--invariant. Therefore all the classes in $E_2$
survive. Moreover, all of the generators of the two factors in the
decomposition \eqref{E2} lift canonically to the cohomology
$H^\bullet(\g[[t]],\V)$, and so we obtain the desired statement.
\end{proof}

\section{Cohomology of the classical Verma module}    
\label{class verma}

Let $\bb$ be a Borel subalgebra of $\g$. Denote by $\tb$ the Lie
subalgebra of $\g[[t]]$ which is the preimage of $\bb$ under the
evaluation homomorphism $\g[[t]] \to \g$. Thus, we can write
$$
\tb = \bb \otimes 1 \bigoplus \g \otimes t \C[[t]].
$$
The Lie group of $\tb$ is the preimage of $B$, the Borel subgroup of
$G$ corresponding to $\bb$ under the evaluation homomorphism $G[[t]]
\to G$. We call it the Iwahori subgroup and denote it by $\wt{B}$.

Decompose $\bb = \h \oplus \n$ into a direct sum of a Cartan
subalgebra and a nilpotent subalgebra. Then we have an analogous
decomposition of $\tb$: $\tb = \h \oplus \tn$, where
\begin{equation}    \label{tn}
\tn = \n \otimes 1 \bigoplus \g \otimes t \C[[t]].
\end{equation}
Given $\la \in \h^*$, let $\C_\la$ be the one-dimensional
representation of the Lie subalgebra $\tb \oplus \C K$ of
$\ghat_{\ka}$, on which $\tn$ acts by $0$, $\h$ acts according to
$\la$ and $K$ acts as the identity. Then the Verma module
$\M_{\la,\ka}$ of highest weight $\la$ at level $\ka$ is by definition
the induced module
$$
\M_{\la,\ka} = \on{Ind}_{\tb \oplus \C K}^{\ghat_{\ka}} \C_{\la}.
$$

Let $H$ be the Cartan subgroup of $G$ whose Lie algebra is $\h$. We
have the categories $HC(\ghat_\ka,H)$ and $HC(\ghat_\ka,\wt{B})$ of
continuous Harish-Chandra modules corresponding to the pairs
$(\ghat_\ka,H)$ and $(\ghat_\ka,\wt{B})$, respectively.

We wish to compute the algebras of self-extensions
$$
\on{Ext}^\bullet_{HC(\ghat_\ka,H)}(\M_{\la,\ka},\M_{\la,\ka}) \qquad
\on{and} \qquad
\on{Ext}^\bullet_{HC(\ghat_\ka,\wt{B})}(\M_{\la,\ka},\M_{\la,\ka}).
$$
Consider the relative cohomology
$$
H^\bullet(\g\ppart,\h,\on{End} \M_{\la,\ka})
$$
with its natural algebra structure. In the same way as in \secref{vac
  module}, we obtain the following result:

\begin{prop}
There are natural isomorphisms of algebras
$$
\on{Ext}^\bullet_{HC(\ghat_\ka,H)}(\M_{\la,\ka},\M_{\la,\ka}) \simeq
\on{Ext}^\bullet_{HC(\ghat_\ka,\wt{B})}(\M_{\la,\ka},\M_{\la,\ka})
\simeq H^\bullet(\g\ppart,\h,\on{End} \M_{\la,\ka}).
$$
\end{prop}

Next, consider the relative cohomology
$$
H^\bullet(\tb,\h,\M_{\la,\ka} \otimes \C_{-\la}) =
H^\bullet(\tn,M_{\la,\ka})_{\la},
$$
where in the right hand side we consider the $\la$--component of the
cohomology $H^\bullet(\tn,\M_{\la,\ka_c})$ with respect to the natural action
of the Cartan subalgebra $\h$. Using Shapiro's lemma as in
\lemref{isom with exts}, we obtain an isomorphism
$$
H^\bullet(\tb,\h,\M_{\la,\ka} \otimes \C_{-\la}) \simeq H^\bullet(\g\ppart,
\h,\on{End} \M_{\la,\ka}).
$$
This isomorphism gives an algebra structure to
$H^\bullet(\tb,\h,\M_{\la,\ka} \otimes \C_{-\la})$.

Thus, to compute the algebras of self-extensions of $\M_{\la,\ka}$
introduced above, we need to compute the relative Lie algebra
cohomologies $H^\bullet(\tb,\h,\M_{\la,\ka} \otimes \C_{-\la})$. The
computation will proceed in the same way as in the case of the vacuum
module. Namely, we will use a spectral sequence whose first term is
the cohomology of the graded version of $\M_{\la,\ka}$. In the rest of
this section we describe the latter.

The $\tb$--module module $\M_{\la,\ka} \otimes \C_{-\la}$ carries a PBW
filtration and the associate graded is isomorphic, as a $\tb$--module,
to
$$
\M^{\on{cl}} = \on{Sym}(\g\ppart/\tb).
$$
We consider the relative cohomology
$H^\bullet(\tb,\h,\M^{\on{cl}})$. This cohomology is computed by the
standard Chevalley complex of Lie algebra cohomology
$$
C^\bullet(\tb,\h,\M^{\on{cl}}) = \left( \M^{\on{cl}} \otimes
\bigwedge{}^\bullet \; \tn^* \right)^{\h} = \left(
\on{Sym}(\g\ppart/\tb) \otimes \bigwedge{}^\bullet \; \tn^*
\right)^{\h}.
$$
The algebra structure on $\M^{\on{cl}}$ gives rise to a graded algebra
structure on the cohomology of this complex. We compute this
cohomology using \thmref{classical}.

Recall the space $C_{\g^*}$ from \secref{qc case} and define
$C^{\on{RS}}_{\g^*,\om}$ as
$$
C^{\on{RS}}_{\g^*,\om} = \Gamma\left(D,\left(\om(0)
\underset{\C^\times}\times C_{\g^*}\right)(-0)\right) \simeq t \left(
t^{-1} \C[[t]] dt \underset{\C^\times}\times C_{\g^*} \right),
$$
where $0$ is the closed point of the disc $D = \on{Spec}\C[[t]]$ and
$(-0)$ indicates sections sending the closed point of $D$ to the
origin of the cone. Using the generators $P_i, i=1,\ldots,\ell$ of
$(\on{Fun} \g^*)^G$ introduced in \secref{qc case}, we obtain an
identification
$$
\on{Fun}(C^{\on{RS}}_{\g^*,\om}) = \C[P_{i,n}]_{i=1,\ldots,\ell;n \geq
  -d_i},
$$
where the $P_{i,n}$'s are the functions on $\g^* \otimes t^{-1}
\C[[t]] dt$ defined by the formula \eqref{Pin}. Thus, we obtain a
coordinate-independent isomorphism
$$
C^{\on{RS}}_{\g^*,\om} \simeq \bigoplus_{i=1}^\ell
\Gamma(D,\om^{\otimes(d_i+1)}(d_i \cdot 0)) = t^{-d_i} \C[[t]]
(dt)^{\otimes(d_i+1)}.
$$

Consider the restrictions of the functions $P_{i,n}$ to
$$
(\g\ppart/\tb)^* \subset (\g\ppart/\g \otimes t\C[[t]])^* \simeq \g^*
\otimes t^{-1} \C[[t]] dt.
$$
Since by construction these functions are $\tb$--invariant, we obtain
a map
$$
\on{Fun}(C^{\on{RS}}_{\g^*,\om}) \to H^0(\tb,\h,\M^{\on{cl}}).
$$

Next, we use the map $\varphi^{\on{cl}}_{\ka_0}$ from \secref{qc case}
to construct a map
$$
(C^{\on{RS}}_{\g^*,\om})^* \to H^1(\tb,\h,\M^{\on{cl}}).
$$

\begin{prop}    \label{classical verma}
There is a canonical isomorphism of graded algebras
$$
H^\bullet(\tb,\h,\M^{\on{cl}}) \simeq
\Omega^\bullet(C^{\on{RS}}_{\g^*,\om}) =
\on{Fun}(C^{\on{RS}}_{\g^*,\om}) \otimes
\bigwedge{}^\bullet(C^{\on{RS}}_{\g,\om})^*.
$$
The right hand side is a free skew-commutative algebra with the even
generators $P_{i,n} \in H^0(\tb,\h,\M^{\on{cl}})$ defined by formula
\eqref{Pin} and the odd generators
$\varphi^{\on{cl}}_{\kappa_0}(P_{i,n}) \in$ $H^1(\tb,\h,\M^{\on{cl}})$
defined by formula \eqref{varphi}, where $i=1,\ldots,\ell$ and $n\geq
-d_i$.
\end{prop}

\begin{proof}
We use the same argument as in the proof of Theorem 1.13 of
\cite{FGT}, where a closely related cohomology,
$H^\bullet(\tb,\h,\on{Sym}(\g\ppart/\tn))$, was computed.

Using the van Est spectral sequence, we obtain that
$$
H^\bullet(\tb,\h,\M^{\on{cl}}) \simeq
H^\bullet_{\wt{B}}(\M^{\on{cl}}), \qquad
H^\bullet(\g[[t]],\g,\V^{\on{cl}}) \simeq
H^\bullet_{G[[t]]}(\V^{\on{cl}}).
$$

By Shapiro's lemma, we have a spectral sequence converging to
$H^\bullet_{\wt{B}}(\M^{\on{cl}})$ whose second term consists of the
cohomologies
$$
H^p_{G[[t]]}(R^q \on{Ind}_{\wt{B}}^{G[[t]]} \M^{\on{cl}}).
$$
Using the residue pairing and a non-degenerate inner
product $\ka_0$ on $\g$ we have an isomorphism
$$
\g\ppart/\tb \simeq \tn \frac{dt}{t} \simeq \tn
$$
as a $\tb$--module, and hence
$$
\M^{\on{cl}} = \on{Sym}(\g\ppart/\tb) \simeq \on{Fun}(\tn).
$$

Therefore we have
$$
R^q \on{Ind}_{\wt{B}}^{G[[t]]} \M^{\on{cl}} = H^q\left(G[[t]]
\underset{\wt{B}}\times \tn,\OO\right),
$$
where $\OO$ is the structure sheaf. But the splitting $G[[t]] = G
\times G^{(1)}$, where $G^{(1)}$ is the first congruence subgroup,
gives rise to the following isomorphism of $G$--equivariant
vector bundles over $G[[t]]/\wt{B} \simeq G/B$:
$$
G[[t]] \underset{\wt{B}}\times \tn \simeq G \underset{{B}}\times \tn
\simeq T^*(G/B) \times (\g \otimes t\C[[t]]),
$$
where the second isomorphism is due to the fact that
$$
G \underset{{B}}\times \n \simeq T^*(G/B),
$$
and the second summand in the direct sum decomposition \eqref{tn} is a
$G$--module giving rise to a trivial vector bundle on
$G/B$.

It follows from the results of \cite{Hess} that
$$
H^q(T^*(G/B),\OO) = 0, \qquad q>0,
$$
and
$$
H^0(T^*(G/B),\OO) = \on{Fun}(\Nil),
$$
where $\Nil \subset \g$ is the nilpotent cone and the functions on
$\Nil$ are pulled back to $T^*(G/B)$ via the moment map $T^*(G/B) \to
\g^* \overset{\ka_0}\longrightarrow \g$ (its image belongs to
$\Nil$). This implies that $R^q \on{Ind}_{\wt{B}}^{G[[t]]}
\M^{\on{cl}} = 0$ for $q>0$. Consider the morphism
\begin{align*}
G[[t]] \underset{\wt{B}}\times \tn &\to T^*(G/B) \underset{\g}\times
\g[[t]], \\
(g,x) &\mapsto ((\ol{g},\ol{g(x)}),g(x)),
\end{align*}
where $\ol{g}$ is the projection of $g \in G[[t]]$ onto $G[[t]]/\wt{B}
\simeq G/B$ and $\ol{g(x)} = g(x) \; \on{mod} \; \g \otimes t\C[[t]]$, so
that $\ol{g(x)} \in T^*_{\ol{g}}(G/B)$. It is clear that this is an
isomorphism, and so we find that
$$
\on{Ind}_{\wt{B}}^{G[[t]]} \M^{\on{cl}} \simeq \on{Fun}\left(G[[t]]
\underset{\wt{B}}\times \tn \right) \simeq \on{Fun} \left(\Nil
\underset{\g}\times \g[[t]]\right).
$$
Therefore we obtain that
$$
H^p_{\wt{B}}(\M^{\on{cl}}) = H^p_{G[[t]]}\left(\on{Fun} \left(\Nil
\underset{\g}\times \g[[t]]\right) \right),
$$
Recall that $\Nil$ is a complete intersection whose ideal in $\on{Fun}
\g$ is the augmentation ideal of the ring of invariant functions
$(\on{Fun} \g)^G$. Thus, the ring of functions on the fiber product
appearing on the right hand side may be resolved by the Koszul complex
$$
\on{Fun} \g \otimes \bigwedge{}^\bullet(P_i)_{i=1,\ldots,\ell}.
$$
Hence we find that $H^\bullet(\tb,\h,\M^{\on{cl}})$ is equal to the
cohomology of the double complex
$$
H^\bullet\left(\g[[t]],\g,\on{Fun}(\g[[t]]) \right) \otimes
\bigwedge{}^\bullet(P_i)_{i=1,\ldots,\ell}.
$$
But we already know the first factor from \thmref{classical}. It is
clear from the description of the generators of this cohomology that
the even generators $P_{i,0}$ get eliminated by the exterior algebra
$\bigwedge{}^\bullet(P_i)_{i=1,\ldots,\ell}$. Hence we obtain that
$H^\bullet(\tb,\h,\M^{\on{cl}})$ is a free skew-commutative algebra
generated by all remaining generators as described in
\thmref{classical}. This gives us the statement of the proposition,
with the shifting of the labeling of the generators due to the fact
that we have used the identification $\tn dt/t \simeq \tn$ in the
course of the proof.
\end{proof}

\section{Self-extensions of the Verma modules}    \label{q verma}

We can now compute the algebras of self-extensions of the Verma
modules $\M_{\la,\ka}$ in the categories of Harish-Chandra modules
introduced above. First, we describe some deformations of the space
$C^{\on{RS}}_{\g^*,\om}$.

Following \cite{BD}, \S~3.8.8, define a $^L\g$--oper with regular
singularity as the equivalence class operators of the form
$$
\pa_t + \dfrac{1}{t}\left(p_{-1} + {\mb v}(t) \right), \qquad {\mb
  v}(t) \in {}^L \bb[[t]]
$$
with respect to the gauge action of $^L N[[t]]$. Define the residue of
this oper as $p_{-1} + {\mb v}(0)$. Clearly, under gauge
transformations by an element $x(t)$ of $^L N[[t]]$ the residue gets
conjugated by $x(0) \in {}^L N$. Therefore its projection onto
$^L\g/{}^L G = {}^L\h/W = \on{Spec} (\on{Fun} {}^L\h)^W$ is
well-defined. Hence the residue is a point in $^L\h/W =
\h^*/W$. Denote the space of $^L\g$--opers with regular singularity by
$\on{Op}^{\on{RS}}_{^L\g}(D)$ and for $\la \in \h^*$ denote by
$\on{Op}^{\on{RS}}_{^L\g}(D)_{\la}$ its subspace of opers with residue
equal to the projection of $\la$ onto $\h^*/W$.

The natural map
$\on{Op}^{\on{RS}}_{^L\g}(D)_{\la} \to \on{Op}_{^L\g}(D^\times)$ is
an embedding and the canonical form of an oper in
$\on{Op}^{\on{RS}}_{^L\g}(D)_{\la}$ is given by the formula
$$
\pa_t + p_{-1} + \sum_{i=1}^\ell t^{-d_i-1} c_j(t) p_j,
\qquad c_i(t) \in \C[[t]].
$$
Moreover, the values $c_i(0)$ are uniquely determined by the
requirement that $\la \in {}^L \h$ is $^L G$--conjugate to
$$
p_{-1} +  \left( c_1(0) + \frac{1}{4} \right) p_1 + \sum_{i=2}^\ell
c_j(0) p_j
$$
(see \cite{BD}, Prop. 3.8.9). It follows that the algebra
$\on{Fun}(\on{Op}^{\on{RS}}_{^L \g}(D))$ has a natural filtration, and
the associated graded algebra is isomorphic to
$\on{Fun}(C^{\on{RS}}_{\g^*,\om})$.

According to the results \cite{F:wak} (see Theorem 12.4, Lemma 9.4 and
Prop. 12.8), we have the following analogue of \thmref{center} for
Verma modules.

\begin{thm}    \label{end of verma}
For any $\la \in \h^*$ there is a canonical isomorphism of algebras
$$
\on{End}_{\ghat_{\ka_c}} \M_{\la,\ka_c} = H^0(\wt\bb,\h,\M_{\la,\ka_c} \otimes
\C_{-\la})\simeq \on{Fun}(\on{Op}^{\on{RS}}_{^L\g}(D)_{-\la-\rho}).
$$
\end{thm}

The main step in the proof is proving that the associated graded of
$(\M_{\la,\ka_c} \otimes \C_{-\la})^{\wt{\bb}}$ with respect to the PBW
filtration is equal to $(\M^{\on{cl}})^{\wt{\bb}}$. We prove this by
constructing $\wt\bb$--invariant vectors in $\M_{\la,\ka_c} \otimes \C_{-\la}$
using central elements of the completed universal enveloping algebra
of $\ghat_{\ka_c}$. Since the center is ``large'', namely, it is
isomorphic to the algebra of functions on $\on{Op}_{^L\g}(D^\times)$,
we construct sufficiently many invariant vectors this way. But we know
from \propref{classical verma} that
$$
(\M^{\on{cl}})^{\wt{\bb}} \simeq \on{Fun}(C^{\on{RS}}_{\g^*,\om}).
$$
Therefore $\on{Spec} \on{End}_{\ghat_{\ka_c}} \M_{\la,\ka_c}$ is an
affine subspace of $\on{Op}_{^L\g}(D^\times)$ on which
$C^{\on{RS}}_{\g^*,\om}$ acts simply transitively. The fact that this
subspace is equal to $\on{Op}^{\on{RS}}_{^L\g}(D)_{-\la-\rho}$ follows
from our knowledge of how the degree zero part of the center acts on
the highest weight vector of $\M_{\la,\ka_c}$ (see \cite{F:wak},
Prop. 12.8).

As in the case of the vacuum module, we now generalize \thmref{end of
verma} to a complete description of the algebra
$H^\bullet(\wt\bb,\h,\M_{\la,\ka_c} \otimes \C_{-\la})$. It turns out
to be isomorphic to the algebra of differential forms on
$\on{Op}^{\on{RS}}_{^L\g}(D)_{-\la-\rho}$.

\begin{thm}    \label{verma}
For each non-zero invariant inner product $\kappa_0$ on $\g$ there is
an isomorphism of graded algebras
\begin{multline}    \label{isom coh1}
\on{Ext}^\bullet_{HC(\ghat_{\ka_c},\wt{B})}(\M_{\la,\ka_c},
\M_{\la,\ka_c}) = H^\bullet(\wt\bb,\h,\M_{\la,\ka_c} \otimes
\C_{-\la}) \\ \simeq \Omega^\bullet(\on{Op}^{\on{RS}}_{^L
\g}(D)_{-\la-\rho}) = \on{Fun}(\on{Op}^{\on{RS}}_{^L
\g}(D)_{-\la-\rho}) \otimes
\bigwedge{}^\bullet(C^{\on{RS}}_{\g^*,\om})^*,
\end{multline}
where $\Omega^\bullet(\on{Op}^{\on{RS}}_{^L \g}(D)_{-\la-\rho})$ is
the algebra of differential forms on $\on{Op}^{\on{RS}}_{^L
\g}(D)_{-\la-\rho}$.  If we rescale $\ka_0$ by $\lambda$, then this
isomorphism gets rescaled by $\la^{-i}$ on the $i$th cohomology.
\end{thm}

\begin{proof}
We apply verbatim the proof of \thmref{vacuum}. We compute the
cohomology $H^\bullet(\tb,\h,\M_{\la,\ka_c} \otimes \C_{-\la})$ using the
spectral sequence associated to the PBW filtration on
$\M_{\la,\ka_c}$. Its first term is the cohomology
$H^\bullet(\tb,\h,\M^{\on{cl}})$ computed in \propref{classical
verma}. We find from \thmref{end of verma} that the zeroth cohomology
part of the first term survives. Next, we construct cohomology classes
in $H^1(\tb,\h,\M_{\la,\ka_c} \otimes \C_{-\la})$ in the same way as in
\thmref{vacuum}, by deforming the module $\M_{\la,\ka_c}$ away from
the critical level.

In the same way as in the case of the vacuum module we construct maps
$$
\wt{\varphi}_{\ka_0}^i: H^i(\tb,\h,\M_{\la,\ka_c} \otimes \C_{-\la}) \to
H^{i+1}(\tb,\h,\M_{\la,\ka_c} \otimes \C_{-\la}).
$$
Considering the symbols of the classes in $H^1(\tb,\h,\M_{\la,\ka_c} \otimes
\C_{-\la})$ obtained by applying $\wt{\varphi}_{\ka_0}^0$ to the
generators of $H^0(\tb,\h,\M_{\la,\ka_c} \otimes \C_{-\la})$, we find that
these are the generators of $H^1(\tb,\h,\M^{\on{cl}})$. Hence these
classes survive.

We then use the product structure to produce cohomology classes of
higher degrees. We check in the same way as in the proof of
\thmref{vacuum} (the second argument) that the product on
$H^\bullet(\tb,\h,\M_{\la,\ka_c} \otimes \C_{-\la})$ is compatible
with that on $H^\bullet(\tb,\h,\M^{\on{cl}})$. Therefore we find that
all of these classes survive. Finally, we show that the generators of
$H^1(\tb,\h,\M_{\la,\ka_c} \otimes \C_{-\la})$ anti-commute using the
same computation as in \remref{vacuum}. This completes the proof.
\end{proof}

Finally, we describe the algebras of self-extensions of the Verma
modules $\M_{\la,\ka}$ of an arbitrary level $\ka \neq \ka_c$.

\begin{prop}
For $\ka \neq \ka_0$ and any $\la \in \h^*$ we have
$$\on{Hom}_{\g_\ka}(\M_{\la,\ka},\M_{\la,\ka}) =
H^0(\tb,\h,\M_{\la,\ka} \otimes \C_{-\la}) = \C$$ and
$$
\on{Ext}^i_{HC(\ghat_{\ka},\wt{B})}(\M_{\la,\ka},
\M_{\la,\ka}) = H^i(\tb,\h,\M_{\la,\ka} \otimes \C_{-\la}) = 0
$$
for all $i>0$.
\end{prop}

\begin{proof}
The proof is the verbatim repetition of the proof of \propref{gen
ka}. First of all, for generic $\ka$ and $\la$ the module
$\M_{\la,\ka}$ is irreducible, according to the results of
\cite{KK}. Therefore it is $\tn$--cofree, hence the result.

For an arbitrary $\ka \neq \ka_c$ and $\la$ we use the same spectral
sequence as in the proof of \thmref{verma} induced by the PBW
filtration on $\M_{\la,\ka}$. Then the first term of this spectral
sequence is the same as for $\ka=\ka_c$, i.e.,
$H^\bullet(\tb,\h,\M^{\on{cl}})$. Next, we find that the first
differential of this spectral sequence differs from the first
differential of the corresponding spectral sequence at $\ka=\ka_c$ by
$h\on{gr} \wt{\varphi}^i_{\ka_0}$, where $\ka = \ka_c + h \ka_0$. But
we know from the proof of \thmref{verma} that the first differential
of the spectral sequence at $\ka=\ka_c$ is equal to $0$, and we know
that $\wt{\varphi}^i_{\ka_0}$ is just the de Rham differential on
$H^\bullet(\tb,\h,\M^{\on{cl}}) \simeq
\Omega^\bullet(C^{\on{RS}}_{\g^*,\om})$. Hence we obtain the desired
assertion.
\end{proof}

\end{document}